\PassOptionsToPackage{unicode}{hyperref}
\PassOptionsToPackage{hyphens}{url}
\PassOptionsToPackage{dvipsnames,svgnames,x11names}{xcolor}
\documentclass[
  12pt]{article}

\usepackage{amsmath,amssymb}
\usepackage{iftex}
\ifPDFTeX
  \usepackage[T1]{fontenc}
  \usepackage[utf8]{inputenc}
  \usepackage{textcomp} 
\else 
  \usepackage{unicode-math}
  \defaultfontfeatures{Scale=MatchLowercase}
  \defaultfontfeatures[\rmfamily]{Ligatures=TeX,Scale=1}
\fi
\usepackage{lmodern}
\usepackage{breakurl}
\usepackage[linesnumbered,boxed,ruled,commentsnumbered]{algorithm2e}
\newtheorem{definition}{Definition}
\newtheorem{theorem}{Theorem}
\newtheorem{lemma}{Lemma}

\newtheorem{proposition}{Proposition}
\newtheorem{remark}{Remark}
\newtheorem{example}{Example}
\usepackage{enumerate}
\ifPDFTeX\else  
\fi
\IfFileExists{upquote.sty}{\usepackage{upquote}}{}
\IfFileExists{microtype.sty}{
  \usepackage[]{microtype}
  \UseMicrotypeSet[protrusion]{basicmath} 
}{}
\makeatletter
\@ifundefined{KOMAClassName}{
  \IfFileExists{parskip.sty}{%
    \usepackage{parskip}
  }{
    \setlength{\parindent}{0pt}
    \setlength{\parskip}{6pt plus 2pt minus 1pt}}
}{
  \KOMAoptions{parskip=half}}
\makeatother
\usepackage{xcolor}
\makeatletter
\ifx\paragraph\undefined\else
  \let\oldparagraph\paragraph
  \renewcommand{\paragraph}{
    \@ifstar
      \xxxParagraphStar
      \xxxParagraphNoStar
  }
  \newcommand{\xxxParagraphStar}[1]{\oldparagraph*{#1}\mbox{}}
  \newcommand{\xxxParagraphNoStar}[1]{\oldparagraph{#1}\mbox{}}
\fi
\ifx\subparagraph\undefined\else
  \let\oldsubparagraph\subparagraph
  \renewcommand{\subparagraph}{
    \@ifstar
      \xxxSubParagraphStar
      \xxxSubParagraphNoStar
  }
  \newcommand{\xxxSubParagraphStar}[1]{\oldsubparagraph*{#1}\mbox{}}
  \newcommand{\xxxSubParagraphNoStar}[1]{\oldsubparagraph{#1}\mbox{}}
\fi
\makeatother

\usepackage{longtable,booktabs,array}
\usepackage{calc} 
\usepackage{etoolbox}
\makeatletter
\patchcmd\longtable{\par}{\if@noskipsec\mbox{}\fi\par}{}{}
\makeatother
\IfFileExists{footnotehyper.sty}{\usepackage{footnotehyper}}{\usepackage{footnote}}
\makesavenoteenv{longtable}
\usepackage{graphicx}
\makeatletter
\def\maxwidth{\ifdim\Gin@nat@width>\linewidth\linewidth\else\Gin@nat@width\fi}
\def\maxheight{\ifdim\Gin@nat@height>\textheight\textheight\else\Gin@nat@height\fi}
\makeatother
\setkeys{Gin}{width=\maxwidth,height=\maxheight,keepaspectratio}
\makeatletter
\def\fps@figure{htbp}
\makeatother

\makeatletter
\@ifpackageloaded{caption}{}{\usepackage{caption}}
\AtBeginDocument{%
\ifdefined\contentsname
  \renewcommand*\contentsname{Table of contents}
\else
  \newcommand\contentsname{Table of contents}
\fi
\ifdefined\listfigurename
  \renewcommand*\listfigurename{List of Figures}
\else
  \newcommand\listfigurename{List of Figures}
\fi
\ifdefined\listtablename
  \renewcommand*\listtablename{List of Tables}
\else
  \newcommand\listtablename{List of Tables}
\fi
\ifdefined\figurename
  \renewcommand*\figurename{Figure}
\else
  \newcommand\figurename{Figure}
\fi
\ifdefined\tablename
  \renewcommand*\tablename{Table}
\else
  \newcommand\tablename{Table}
\fi
}
\@ifpackageloaded{float}{}{\usepackage{float}}
\floatstyle{ruled}
\@ifundefined{c@chapter}{\newfloat{codelisting}{h}{lop}}{\newfloat{codelisting}{h}{lop}[chapter]}
\floatname{codelisting}{Listing}

\makeatother
\makeatletter
\@ifpackageloaded{caption}{}{\usepackage{caption}}
\@ifpackageloaded{subcaption}{}{\usepackage{subcaption}}
\makeatother

\ifLuaTeX
  \usepackage{selnolig}  
\fi
\usepackage[]{natbib}
\usepackage{bookmark}

\IfFileExists{xurl.sty}{\usepackage{xurl}}{} 
\hypersetup{
  pdftitle={Title},
  pdfauthor={Author 1; Author 2},
  pdfkeywords={3 to 6 keywords, that do not appear in the title},
  colorlinks=true,
  linkcolor={blue},
  filecolor={Maroon},
  citecolor={Blue},
  urlcolor={Blue},
  pdfcreator={LaTeX via pandoc}}

\newcommand{\anon}{1}

\begin{document}

\def\spacingset#1{\renewcommand{\baselinestretch}%
{#1}\small\normalsize} \spacingset{1}


\if1\anon
{
  \title{\bf Sliced Space-filling Design with Mixtures}
  \author{Zikang Xiong\quad  Hong Qin\thanks{
    The authors gratefully acknowledge the National Natural Science Foundation of China (No.12371261) and the Fundamental Research Funds for the Central Universities, Zhongnan University of Economics and Law (2722024BQ064).\\
    Corresponding author: Hong Qin (\texttt{qinhong@ccnu.edu.cn})}\quad Yuning Huang\hspace{.2cm} \\
    School of Statistics and Mathematics, \\
    Zhongnan University of Economics and Law, Wuhan, China, 430073\\
    and \\
    Jianhui Ning \\
    School of Mathematics and Statistics,\\ Central China Normal University, Wuhan, China, 430079}
  \maketitle
} \fi

\if0\anon
{
  \bigskip
  \bigskip
  \bigskip
  \begin{center}
    {\LARGE\bf Title}
\end{center}
  \medskip
} \fi

\bigskip
\begin{abstract}
In this paper, we proposes the construction methods of sliced space-filling design when the quantitative factors are mixture components. Leveraging the representative points framework for distribution and energy distance decomposition theory, this paper proposes three methods for constructing sliced representative points and establishes their distributional convergence. Furthermore, one-shot and sequential algorithms for generating sliced space-filling mixture design for experiments with process variables are presented with convergence proofs. Compared to existing methods, the proposed sliced space-filling mixture design exhibits greater flexibility in subdesign run sizes and broader applicability to constrained experimental regions. Moreover, numerical results confirm its marked advantages in both space-filling performance and predictive accuracy.
\end{abstract}

\noindent%

\vfill

\newpage
\spacingset{1.8} 

\section{Introduction}
\label{sec1}
Many mixture experiments not only involve variables of mixture components, 
but also take into account process variables additionally \citep{Cornell2002, Anderson}, such experimental platforms, temperature, pressure. In practice, experimenters focus more on the response surface profiles of mixture variables across different combinations of process variable levels, which is a scenario that inherently involves both quantitative and qualitative variables.

Sliced space-filling designs are widely adopted as model-robust experimental frameworks when both quantitative and qualitative factors coexist. Current literature \citep{sliced2009,slhd2012,uslhd2016,SLHDr,MaxPror} primarily addresses their construction for quantitative variables within hypercube. While inverse transformation methods \citep{Wang1990,Fang.2000} could theoretically extend these designs to mixture experimental region, they suffer from critical limitations: (1) the degradation of space-filling property, and (2) the difficulty in solving inverse transformation functions for mixture experiments with complex constraints.
To enable direct construction within mixture regions, \cite{mixslice} developed an algorithm for sliced designs on the standard simplex using sphere packing \citep{spdm,He2017} and coset decomposition \citep{Heslice}. Alternatively, \cite{FFF2018} proposed approximate constructions for irregular regions via hierarchical clustering and weighted maximum projection criteria. Both approaches exhibit limited dimensionality scalability and inflexible allocation of points in the subdesigns.

This paper establishes a decomposition theory for energy distance within the representative points framework \citep{Fang2022RP, Mak2018}, which characterizes the relationship between the representativeness of a full point set and its partitioned subsets. To simultaneously ensure representativeness at both levels, we propose three methods for constructing representative point set with sliced structure and the resulting full point sets and subsets provably converge to the target distribution. Furthermore, within the hybrid energy distance criterion and Majorization-Minimization algorithm framework, we develop one-shot and sequential construction algorithms for sliced space-filling designs with mixtures in convex constrained regions, establishing their convergence properties. The representative points approach offers three key advantages:
(1) avoidance of transformation functions required by inverse transform methods, facilitating design construction under complex constraints;
(2) free selection of the number of experimental points in each subdesign;
(3) superior uniformity in both full design and subdesigns.

The subsequent sections of this paper are organized as follows. Section 2 introduces basic concepts related to mixture experiments with process variables and sliced space-filling design.   Section 3 proposes some construction theories and algorithms for representative points with sliced structure under the energy distance criterion and applies them to the construction problem of sliced space-filling designs with mixtures. Section 4 compares the uniformity of the sliced space-filling designs with mixtures constructed by different methods and their predictive performance in modeling. The final section provides conclusions and future research directions. All the technical proofs are in the supplementary material.

\section{Preliminary}

\subsection{Mixture experiments with process variables}\label{sec2.1}

\begin{definition}\label{Defi1}
Suppose that there are $p$ mixture components $x_1, \ldots, x_p$ and $q$  process variables $z_1, \ldots, z_q$ in a mixture system. The following set
\begin{equation*}
\mathcal{X}_{\mathcal{Z}} = \left\{  (\boldsymbol{x}^T, \boldsymbol{z}^T)^T \mid \boldsymbol{x} = (x_1, x_2, \cdots, x_p)^{\mathrm{T}} \in \mathcal{X} \subseteq T_p, \boldsymbol{z} = (z_1, z_2, \cdots, z_q)^{\mathrm{T}} \in \mathcal{Z} \right\} 
\end{equation*}
is called the mixture experimental region with process variables, where the standard simplex
\begin{equation}\label{standard simplex}
    T_p = \left\{ \boldsymbol{x} = (x_1, \ldots, x_p)^T \mid x_i \geq 0, i = 1, \ldots, p, \sum_{i = 1}^p x_i = 1 \right\} 
\end{equation}
is the unconstrained mixture experimental region. If there is an additional constraint set $S$ for the mixture component $\boldsymbol{x}$, the resulting experimental region becomes $\mathcal{X} = T_p \bigcap S$, then $\mathcal{X}_{\mathcal{Z}}$ is called the constrained mixture experimental region with process variables.
\end{definition}

In practice, a common constraint set $S$ is linear with respect to the mixture components \citep{zhangli.2021,Fang.2018}, such as upper and lower bound constraints, linear equality and inequality constraints\citep{MixofMix}. Since the constraints in the standard simplex $T_p$ are also linear, the mixture region under general linear constraints can be expressed as
\begin{equation}\label{TLP}
   T_{L,p}=\left\{ \boldsymbol{x}=(x_1,\ldots,x_p)^T\mid A\boldsymbol{x}\leq\boldsymbol{b}, C\boldsymbol{x}=\boldsymbol{d}\right\} . 
\end{equation}

\begin{figure}[!htb]
    \centering
    \includegraphics[width=0.6\linewidth]{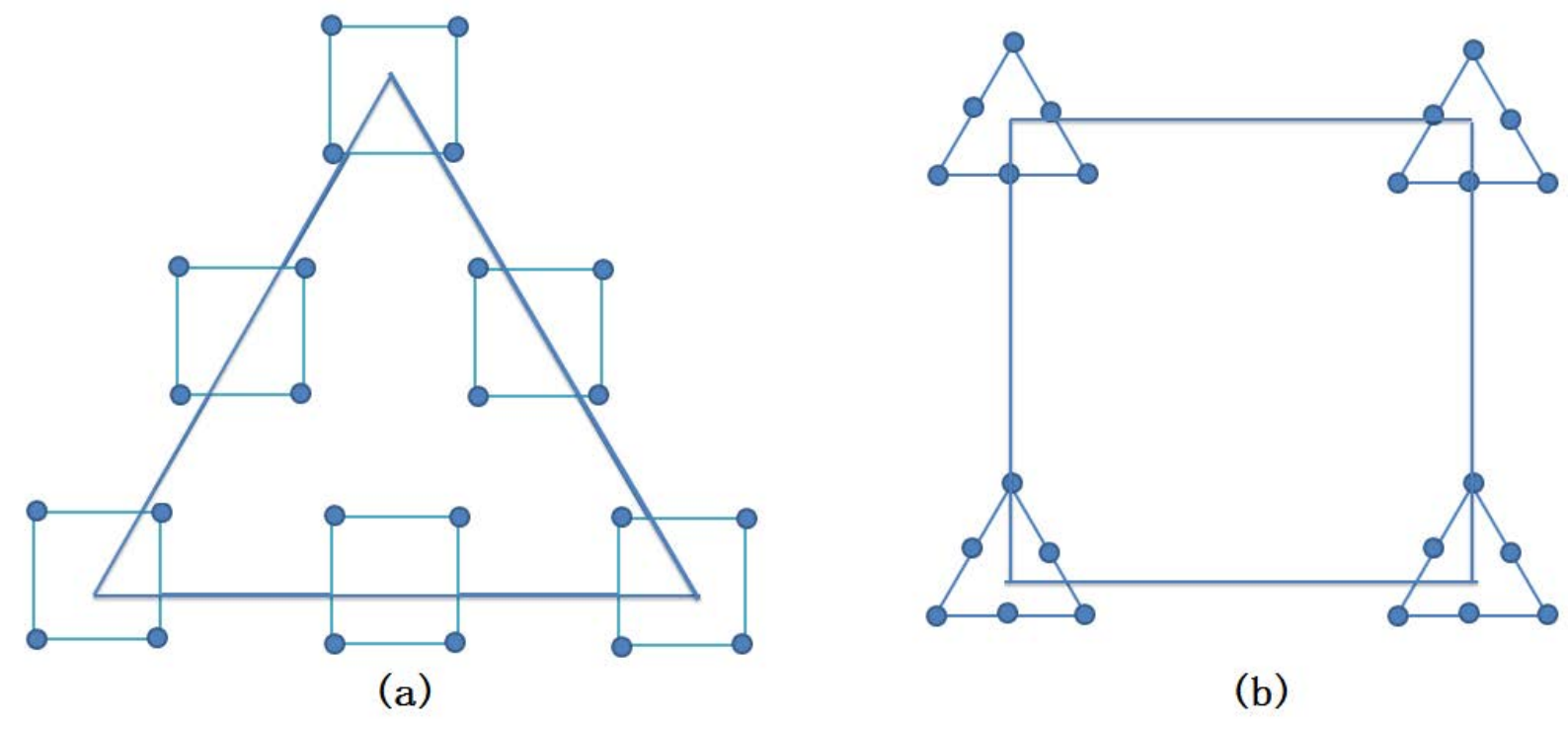}
    \caption{Combined construction scheme for mixture design with process variables ($p=3,\ q=2,\ s_1=s_2=2$).}
    \label{fig1}
\end{figure}

Assume that the number of levels of the process variables are $s_1,\ldots,s_q$ respectively, the experimental region $\mathcal{Z}$ comprises $K = \prod_{i=1}^q s_i$ discrete points formed by all level combinations. The most commonly used mixture design with process variables is constructed by the following combination method: at each experimental point of the mixture design, nest a full (or partial) factorial design of process variables, as shown in Figure\ref{fig1}(a). This design can also be regarded as nesting the same mixture design under each level combination of the factorial design of process variables, as shown in Figure\ref{fig1}(b). 

\subsection{Sliced space-filling design} \label{sec2.2}

Consider an experimental region $\mathcal{X}$ for quantitative variables $\boldsymbol{x}$ and $K$ level combinations of qualitative variables $\boldsymbol{z}$. The quantitative component $\mathcal{P}_n$ of the design $\mathcal{P}=(\mathcal{P}_n, \mathcal{P}_{\boldsymbol{z}})$ partitions according to qualitative variable levels:
\begin{equation}\label{designpartition}
\mathcal{P}_n = \bigcup_{k = 1}^K \mathcal{P}_{n_k}, \quad \sum_{k = 1}^K n_k = n
\end{equation}
where $\mathcal{P}_{n_k}$ denotes the subdesign corresponding to the $k$-th level combination of $\boldsymbol{z}$. The core objective of sliced space-filling designs is to ensure both the full design $\mathcal{P}_n$ and all subdesigns $\mathcal{P}_{n_k} (k = 1, \ldots, K)$ exhibit strong space-filling properties throughout $\mathcal{X}$. 

By introducing the idea of the sliced space-filling design into the mixture experiment, it is required that the mixture design with process variables satisfies that both the full design and the subdesigns have good space-filling properties in the space where the mixture components are located. Such a design is called the sliced space-filling design with mixtures.

\section{Sliced space-filling design with mixtures based on representative points method}\label{sec3}

Uniform designs can be characterized as representative points of the uniform distribution $F_U(\mathcal{X})$ over experimental region $\mathcal{X}$. As discussed in Section \ref{sec2.2} for sliced space-filling designs with mixtures, both the full mixture design $\mathcal{P}_n$ and each subdesign $\mathcal{P}{n_k}$ are required to exhibit space-filling properties throughout $\mathcal{X}$. Thus, full design and all subdesigns are representative points of $F_U(\mathcal{X})$.

\subsection{Energy distance and its decomposition}

\begin{definition}\label{Defi4}
Let $F$ be a continuous distribution function defined in $\mathcal{X}\subseteq \mathbb{R}^p$,  $F_{\mathcal{P}_n}$ be the empirical distribution function of the point set $\mathcal{P}_n = \{ \boldsymbol{x}_i \}_{i = 1}^n \subseteq \mathcal{X}$, and $Y, Y' \stackrel{i.i.d.}{\sim} F$. The energy distance between the distribution functions $F$ and $F_{\mathcal{P}}$ is (assuming these expectations are finite)
\begin{equation}\label{edinms}
  E(F,F_{\mathcal{P}_n})=\frac{2}{n}\sum_{i=1}^n\mathbb{E}\|\boldsymbol{x}_i-Y\|_2-\frac{1}{n^2}
\sum_{i=1}^n\sum_{j=1}^n\|\boldsymbol{x}_i-\boldsymbol{x}_j\|_2-\mathbb{E}\|Y-Y'\|_2.
\end{equation}
\end{definition}

\begin{definition}\label{Defi5}
The point set $\mathcal{P}_n^* = \{ \boldsymbol{x}_i^* \}_{i = 1}^n$ is called the representative points of the distribution function $F$ under the energy distance criterion, if
\begin{align}
 E(F,F_{\mathcal{P}_n^*})&=\underset{\mathcal{P}_n \subseteq \mathcal{X}}{\textup{min}} \ E(F,F_{\mathcal{P}_n})\nonumber\\
&\Leftrightarrow\underset{\mathcal{P}_n \subseteq \mathcal{X}}{\textup{min}} \left\{\frac{2}{n}\sum_{i=1}^n\mathbb{E}\|\boldsymbol{x}_i-Y\|_2-\frac{1}{n^2}
  \sum_{i=1}^n\sum_{j=1}^n\|\boldsymbol{x}_i-\boldsymbol{x}_j\|_2 \label{edrp}\right\}.
\end{align}
\end{definition}

The following Lemma \ref{Lemma1} and \ref{Lemma2} show that the energy distance can be used to measure whether two distribution functions are equal in Euclidean space, which provides a basis for minimizing the energy distance to construct representative points of the distribution function.

\begin{lemma}\label{Lemma1}$\citep{Eddc}$
$ E(F, F_{\mathcal{P}_n}) \geq 0 $, and the equality holds if and only if $ F = F_{\mathcal{P}_n} $.
\end{lemma}

\begin{lemma}\label{Lemma2}$\citep{Mak2018}$
Let $ F_{\mathcal{P}_n^*} $ be the empirical distribution function of the representative point set generated by  $\eqref{edrp}$, then $ F_{\mathcal{P}_n^*} \stackrel{d}{\longrightarrow} F $.
\end{lemma}

To construct the sliced space-filling design with mixtures, we now study the relationship between the energy distance of the full representative point set and the energy distances of the subsets generated by point set partitioning. 

\begin{theorem}\label{Theorem1}
Suppose that the full point set $\mathcal{P}_n = \{ \boldsymbol{x}_i \}_{i=1}^n$ is divided into $K$ subsets $\mathcal{P}_{n_1}, \ldots, \mathcal{P}_{n_K}$ according to  \eqref{designpartition}, then the energy distance corresponding to the full point set can be decomposed as follows:
\begin{equation}\label{decoth}
    \begin{aligned}
E\left(F, F_{\mathcal{P}_n}\right) &= \sum_{k=1}^K \frac{n_k}{n} E\left(F, F_{\mathcal{P}_{n_k}}\right) - \sum_{1 \leq k_1 < k_2 \leq K} \frac{n_{k_1}n_{k_2}}{n^2} E\left(F_{\mathcal{P}_{n_{k_1}}}, F_{\mathcal{P}_{n_{k_2}}}\right) \\
&= \sum_{k=1}^K \frac{n_k}{n} E\left(F, F_{\mathcal{P}_{n_k}}\right) - \sum_{k_1=1}^K \sum_{k_2=1}^K \frac{n_{k_1}n_{k_2}}{2n^2} E\left(F_{\mathcal{P}_{n_{k_1}}}, F_{\mathcal{P}_{n_{k_2}}}\right)\\
&=\sum_{k=1}^K \frac{n_k}{n} \left[E\left(F, F_{\mathcal{P}_{n_k}}\right)-E \left( F_{\mathcal{P}_n}, F_{\mathcal{P}_{n_k}} \right)\right].
\end{aligned}
\end{equation}
In particular, when $K = 2$,
\begin{equation}\label{K2}
   E\left(F, F_{\mathcal{P}_n}\right) = \frac{n_1}{n} E\left(F, F_{\mathcal{P}_{n_1}}\right) + \frac{n_2}{n} E\left(F, F_{\mathcal{P}_{n_2}}\right) - \frac{n_1n_2}{n^2} E\left(F_{\mathcal{P}_{n_1}}, F_{\mathcal{P}_{n_2}}\right). 
\end{equation}
\end{theorem}

\subsection{Construction methods of representative point sets with sliced structures}
\label{sec3.2}
Energy distance serves as a critical metric for evaluating how well a point set represents the target distribution function. A smaller value indicates stronger representativeness. The energy distance decomposition theory in Theorem \ref{Theorem1} provides the basis for analyzing intrinsic representativeness relationships between full point sets and their subsets. We further propose three methods for constructing representative point sets with sliced structure.

\subsubsection{Construct the full point set based on subsets}\label{method1}

Theorem \ref{Theorem1} concludes:
To make the full point set $\mathcal{P}_n$ well represent target distribution $F$, two conditions suffice: (i) enhance representativeness of every subset $\mathcal{P}_{n_1}, \ldots, \mathcal{P}_{n_K}$; (ii) increase pairwise differences between subsets.
Lemma \ref{Lemma1} and Theorem \ref{Theorem1} imply Proposition \ref{Prop1}:
the energy distance of $\mathcal{P}_n$ is dominated by the subsets' energy distances.

\begin{proposition}\label{Prop1}  $0 \leq E(F, F_{\mathcal{P}_n}) \leq \sum_{k=1}^K \frac{n_k}{n} E\left(F, F_{\mathcal{P}_{n_k}}\right).$
\end{proposition}

Let $ n_k(n) $ denote the number of sample points contained in the subset $ \mathcal{P}_{n_k} $, and it is a function of the total sample size $ n $, $\forall\ k = 1, \ldots, K $. Based on Proposition \ref{Prop1}, we can further derive the relationship between the convergence in distribution of the full point set and the subsets to the target distribution.

\begin{theorem}\label{Theorem2}
Suppose $ n = \sum_{k=1}^K n_k(n) $ and $ \lim_{n \to \infty} \frac{n_k(n)}{n} = c_k $ for $ k = 1, \ldots, K $. Let $ S = \{ k \mid c_k > 0, \, k \in \{1, \ldots, K\} \} $. If $\forall\ k \in S$, the subset $ \mathcal{P}_{n_k(n)}^* $ independently generated by  $\eqref{edrp}$, then $  F_{\mathcal{P}_n^*} \stackrel{d}{\longrightarrow} F $.
\end{theorem}

Theorem \ref{Theorem2} demonstrates that for large total sample sizes $n$, the representativeness of the full design is determined by the representativeness of subsets with asymptotically dominant sample sizes. This property proves particularly advantageous when augmenting existing designs \citep{ssp}, as ensuring adequate representativeness in the subsets alone suffices to guarantee full design representativeness without requiring strict distributional homogeneity across subsets. 

\begin{proposition}\label{Prop2}
When $ \mathcal{P}_{n_1} = \cdots = \mathcal{P}_{n_K}=\mathcal{P}_{\frac{n}{K}}$, we have
$E(F, F_{\mathcal{P}_n}) = E\left(F, F_{\mathcal{P}_{\frac{n}{K}}}\right).$
\end{proposition}

Proposition \ref{Prop2} shows that when the subsets are identical optimal representative point sets, the representativeness of the full point set is equivalent to that of a single subset. Therefore, when the sample size is limited, after generating the subsets separately, it is necessary to adjust them so that they not only have good representativeness for the target distribution but also exhibit differences between their distributions. Such post-hoc adjustment is relatively challenging. An alternative approach is to generate subsets incrementally based on  \eqref{K2} using the following sequential algorithm\citep{ssp}:
\begin{enumerate}[$(i)$]
    \item Solve the optimization problem in  \eqref{edrp} to construct the subset $ \mathcal{P}_{n_1}^* $, and set $ \mathcal{P}_{n_c} = \mathcal{P}_{n_1}^* $ and $ n_c = n_1 $;
\item For $ k = 2, \ldots, K $, repeat the following process: Solve the following optimization problem to generate the point set $ \mathcal{P}_{n_k}^* $:  
    \begin{equation}\label{seq}
        \min_{\mathcal{P}_{n_k} \subseteq \mathcal{X}} E\left(F, F_{\mathcal{P}_{n_k}}\right) - \frac{n_c}{n_c + n_k} E\left(F_{\mathcal{P}_{n_c}}, F_{\mathcal{P}_{n_k}}\right),
    \end{equation}  
    Update $ \mathcal{P}_{n_c} = \mathcal{P}_{n_c} \cup \mathcal{P}_{n_k}^* $ and $ n_c = n_c + n_k $.  
\end{enumerate}

The optimization problem in \eqref{seq} reveals that the above process enhances the representativeness of the full point set at the expense of compromising subsequent subsets. Consequently, the representativeness of these later-generated subsets becomes degraded, particularly when substantial differences exist in subset sample sizes.

Theorem \ref{Theorem2} demonstrates that subsets with smaller sample sizes exhibit greater susceptibility to representational compromise, as evidenced by the coefficient $\frac{n_c}{n_c + n_k}$ in \eqref{seq}'s second term. We therefore recommend constructing subsets in ascending sample size order, i.e., sorting such that $n_{(1)} \leq n_{(2)} \leq \cdots \leq n_{(K)}$. When sorting is infeasible, an alternative strategy involves explicitly trading off representativeness between the global point set and subsets during sequential generation, which will be discussed in section \ref{method3}.  

\subsubsection{Construct subsets from the full point set}\label{method2}

The second intuitive approach involves first generating a full representative point set and then partitioning it into representative subsets:
\begin{enumerate}[$(i)$]
    \item Generate the full representative point set $ \mathcal{P}_n^* $ by \eqref{edrp};  
    \item Partition $ \mathcal{P}_n^*$ by solving the following optimization problem:  
    \begin{equation}\label{par0}
    \begin{split}
    \min_{\mathcal{P}_{n_1}, \ldots, \mathcal{P}_{n_K} \subseteq \mathcal{P}_n^*} &\quad \sum_{k=1}^K \frac{n_k}{n} E \left( F, F_{\mathcal{P}_{n_k}} \right) \\
    \text{s.t.} &\quad \bigcup_{k=1}^K \mathcal{P}_{n_k} = \mathcal{P}_n^*, \quad \sum_{k=1}^K n_k = n.
    \end{split}
    \end{equation}
\end{enumerate} 

\begin{theorem}\label{Theorem3}
Suppose $ n = \sum_{k=1}^K n_k(n) $ and $ \lim_{n \to \infty} \frac{n_k(n)}{n} = c_k > 0 $ for all $ k = 1, \ldots, K $. If the point set $ F_{\mathcal{P}_n^*} $ generated by $\eqref{edrp}$, then the subsets generated by $\eqref{par0}$ satisfy $F_{\mathcal{P}_{n_k(n)}^*}  \stackrel{d}{\longrightarrow} F $, $\forall\  k = 1, \ldots, K $.
\end{theorem}

After generating the full point set $ \mathcal{P}_n^* $, the energy distance $ E\left(F, F_{\mathcal{P}_n^*}\right) $ becomes a constant. According to Theorem \ref{Theorem1}, the target distribution $F$ can be removed from the optimization problem \eqref{par0}, which is equivalent to the following optimization problems \eqref{par1} and \eqref{par2}.
\begin{equation}\label{par1}
\begin{split}
\min_{\mathcal{P}_{n_1}, \ldots, \mathcal{P}_{n_K} \subseteq \mathcal{P}_n^*} &\quad \sum_{k_1=1}^K \sum_{k_2=1}^K \frac{n_{k_1}n_{k_2}}{2n^2} E\left(F_{\mathcal{P}_{n_{k_1}}}, F_{\mathcal{P}_{n_{k_2}}}\right) \\
\text{s.t.} &\quad \bigcup_{k=1}^K \mathcal{P}_{n_k} = \mathcal{P}_n^*, \quad \sum_{k=1}^K n_k = n.
\end{split}
\end{equation}  
\begin{equation}\label{par2}
\begin{split}
\min_{\mathcal{P}_{n_1}, \ldots, \mathcal{P}_{n_K} \subseteq \mathcal{P}_n^*} &\quad \sum_{k=1}^K\frac{n_k}{n} E \left( F_{\mathcal{P}_n^*}, F_{\mathcal{P}_{n_k}} \right) \\
\text{s.t.} &\quad \bigcup_{k=1}^K \mathcal{P}_{n_k} = \mathcal{P}_n^*, \quad \sum_{k=1}^K n_k = n. \end{split}
\end{equation}

The aforementioned algorithms can not only generate representative point sets with sliced structure but also be effectively applied to several key problems in the field of machine learning, such as data partitioning (e.g., splitting into training/validation/test sets) and data subsampling (e.g., constructing core sets or alleviating class imbalance). In the case of $ K=2 $, \cite{twinning} and \cite{SPlit} respectively proposed algorithms to partition the complete data into training and test sets by solving the optimization problems in \eqref{par1} and \eqref{par2}. They established that exact solution of this discrete optimization problem is NP-hard and consequently proposed sequential point-selection strategies to partition the full point set into two subsets. This partitioning step can be iterated to obtain any number of subsets K in \eqref{par1} and \eqref{par2}. For algorithmic details, consult \cite{twinning,SPlit} and their corresponding R packages: \texttt{twinning} \citep{twinningr} and \texttt{SPlit} \citep{SPlitr}. For unbalanced subset sample sizes, subsets should be constructed in ascending order of sample size, analogous to the approach described in section \ref{method1}.

\subsubsection{ Construct sliced representative point set with a hybrid criterion}
\label{method3}
Consider the following hybrid criterion that combines the energy distance corresponding to the full point set and the sample proportion-weighted energy distances corresponding to the subsets:  
\begin{equation}
E^\lambda(F, F_{\mathcal{P}_n}) = \lambda E(F, F_{\mathcal{P}_n}) + (1-\lambda) \sum_{k=1}^K \frac{n_k}{n} E\left(F, F_{\mathcal{P}_{n_k}}\right), 
\end{equation}  
where $\lambda \in [0,1]$. From the decomposition of the energy distance in Theorem \ref{Theorem1}, we have:  
\begin{equation*}
E^\lambda(F, F_{\mathcal{P}_n}) = \sum_{k=1}^K \frac{n_k}{n} E\left(F, F_{\mathcal{P}_{n_k}}\right) - \lambda \sum_{k_1=1}^K \sum_{k_2=1}^K \frac{n_{k_1}n_{k_2}}{2n^2} E\left(F_{\mathcal{P}_{n_{k_1}}}, F_{\mathcal{P}_{n_{k_2}}}\right). 
\end{equation*}  
$E^\lambda(F, F_{\mathcal{P}_n})$ depends only on the subsets. Thus, by solving the following optimization problem:  
\begin{equation}\label{mixcri}
\min_{\mathcal{P}_{n_1}, \ldots, \mathcal{P}_{n_K} \subseteq \mathcal{X}} E^\lambda(F, F_{\mathcal{P}_n}), 
\end{equation}  
we can simultaneously obtain $K$ subsets with good representativeness and the full point set formed by their combination.  

\begin{theorem}\label{Theorem4}
Let $ n = \sum_{k=1}^K n_k(n) $ and $ \lim_{n \to \infty} \frac{n_k(n)}{n} = c_k > 0 $ for all $ k = 1, \ldots, K $. Then, the point set $ \mathcal{P}_n^* = \bigcup_{k=1}^K \mathcal{P}_{n_k(n)}^* $ generated based on \eqref{mixcri} satisfies $ F_{\mathcal{P}_{n}^*}  \stackrel{d}{\longrightarrow} F $ and $F_{\mathcal{P}_{n_k(n)}^*}  \stackrel{d}{\longrightarrow} F $, $\forall\ k = 1, \ldots, K $.
\end{theorem}

In the method of section \ref{method1}, the representative point set is generated by sequentially solving the optimization problem in \eqref{seq}. This sequential sampling method offers greater flexibility in practical applications: it allows early termination of the experiment when sufficient evidence is obtained to prove efficacy, inefficacy, or insignificant effects, thereby avoiding unnecessary subsequent data collection.  
As mentioned earlier, the representativeness of sequentially generated subsets and the full point set is affected by the sample sizes of the subsets, and this issue is particularly prominent when the sample sizes are unbalanced. To address this, we further consider adopting a sequential strategy to generate subsets incrementally under the hybrid criterion $ E^\lambda(F, F_{\mathcal{P}_n}) $. Specifically, we replace \eqref{seq} with the following optimization problem:  
\begin{equation}\label{seqmix}
\min_{\mathcal{P}_{n_k} \subseteq \mathcal{X}} E\left(F, F_{\mathcal{P}_{n_k}}\right) - \lambda_k \frac{n_c}{n_c + n_k} E\left(F_{\mathcal{P}_{n_c}}, F_{\mathcal{P}_{n_k}}\right),
\end{equation}  
where the original parameter $ \lambda $ is replaced by $ \lambda_k $, which is a balancing parameter that can be dynamically adjusted with the sequential stage.  

\subsection{Construction algorithms for sliced space-filling design with mixtures}

In this section, we present specific algorithms for sliced space-filling designs with mixtures in $ \mathcal{X} $. For this problem, the target distribution is the uniform distribution $ F_U $ over $ \mathcal{X} $. Since common mixture experimental regions $ \mathcal{X} $ in practice, such as $ T_{L,p} $ in \eqref{TLP}, are bounded convex subsets with zero curvature in Euclidean space, and all moments of the uniform distribution in this region exist. Therefore, representative points methods in section \ref{method2} can be applied, the empirical distributions of the experimental points in both the full design and the subdesigns converge to the uniform distribution $F_U$. 

Since the processes of generating representative points via energy distance optimization in methods of section \ref{method1} and \ref{method2} are special cases of method in section \ref{method3} when $\lambda = 1$ or $\lambda_k = 1$. Therefore, only the detailed algorithms process of the methods based on hybrid energy distance are given in this paper. 

According to \eqref{mixcri} and Definition \ref{Defi5}, we first generate uniform samples $\{\boldsymbol{y}_m\}_{m=1}^N$ in $\mathcal{X}$ to approximate the expectations in the criterion. Then sliced space-filling design with mixtures $\mathcal{P}_n^* = \bigcup_{k=1}^K \mathcal{P}_{n_k}^*$ can be obtained by solving:
\begin{equation}\label{ophy}
\begin{split}
\min_{\mathcal{P}_{n_1}, \ldots, \mathcal{P}_{n_K} \subseteq \mathcal{X}} h(\mathcal{P}_n) &= \frac{2}{nN} \sum_{m=1}^N \sum_{k=1}^K \sum_{i \in \mathcal{I}_{n_k}} \|\boldsymbol{x}_i - \boldsymbol{y}_m\|_2 \\
&\quad - \sum_{k=1}^K \left( \frac{1-\lambda}{nn_k} + \frac{\lambda}{n^2} \right) \sum_{i,j \in \mathcal{I}_{n_k}} \|\boldsymbol{x}_i - \boldsymbol{x}_j\|_2 \\
&\quad - \frac{\lambda}{n^2} \sum_{k_1=1}^K \sum_{\substack{k_2=1 \\ k_2 \neq k_1}}^K \sum_{i \in \mathcal{I}_{n_{k_1}}} \sum_{j \in \mathcal{I}_{n_{k_2}}} \|\boldsymbol{x}_i - \boldsymbol{x}_j\|_2,
\end{split}
\end{equation}  
where $\mathcal{I}_{n_k}$ denotes the index set of points in $\mathcal{P}_{n_k}$ within the full set $\mathcal{P}_n$.

To solve the non-convex optimization problem in \eqref{ophy} iteratively within the Majorization-Minimization (MM) algorithm framework \citep{MMcon}, we need to construct a surrogate function for the objective function $h(\mathcal{P}_n)$ at the current iteration $\mathcal{P}_n^{(t)}$ in each iteration step, which is presented in Theorem \ref{Theorem6} below.

\begin{theorem}\label{Theorem6}
Suppose $\{\boldsymbol{y}_m\}_{m=1}^N$ and $\mathcal{P}_n^{(t)} \subseteq \mathcal{X}$. Define the following function:  
\begin{equation}
\begin{split}
g(\mathcal{P}_n \mid \mathcal{P}_n^{(t)}) &= \frac{1}{nN} \sum_{m=1}^N \sum_{k=1}^K \sum_{i \in \mathcal{I}_{n_k}} \left( \frac{\|\boldsymbol{x}_i - \boldsymbol{y}_m\|_2^2}{\|\boldsymbol{x}_i^{(t)} - \boldsymbol{y}_m\|_2} + \|\boldsymbol{x}_i^{(t)} - \boldsymbol{y}_m\|_2 \right) - \sum_{k=1}^K \left( \frac{1-\lambda}{nn_k} \right. \\
&\quad + \left. \frac{\lambda}{n^2} \right) \sum_{\substack{i,j \in \mathcal{I}_{n_k} \\ i \neq j}} \left[ \|\boldsymbol{x}_i^{(t)} - \boldsymbol{x}_j^{(t)}\|_2 + \frac{2(\boldsymbol{x}_i^{(t)} - \boldsymbol{x}_j^{(t)})^T (\boldsymbol{x}_i - \boldsymbol{x}_i^{(t)})}{\|\boldsymbol{x}_i^{(t)} - \boldsymbol{x}_j^{(t)}\|_2} \right] \\
&\quad - \frac{\lambda}{n^2} \sum_{k_1=1}^K \sum_{\substack{k_2=1 \\ k_2 \neq k_1}}^K \sum_{i \in \mathcal{I}_{n_{k_1}}} \sum_{j \in \mathcal{I}_{n_{k_2}}} \left[ \|\boldsymbol{x}_i^{(t)} - \boldsymbol{x}_j^{(t)}\|_2 + \frac{2(\boldsymbol{x}_i^{(t)} - \boldsymbol{x}_j^{(t)})^T (\boldsymbol{x}_i - \boldsymbol{x}_i^{(t)})}{\|\boldsymbol{x}_i^{(t)} - \boldsymbol{x}_j^{(t)}\|_2} \right].
\end{split}
\end{equation}  
Then, $g(\mathcal{P}_n \mid \mathcal{P}_n^{(t)})$ is a surrogate function of the objective function $h(\mathcal{P}_n)$ in $\eqref{ophy}$ at the current iteration $\mathcal{P}_n^{(t)}$, satisfying:  
$(i)$ $g(\mathcal{P}_n \mid \mathcal{P}_n^{(t)}) \geq h(\mathcal{P}_n)$ for all $\mathcal{P}_n \subseteq \mathcal{X}$;  
$(ii)$ $g(\mathcal{P}_n^{(t)} \mid \mathcal{P}_n^{(t)}) = h(\mathcal{P}_n^{(t)})$.  
Further, let $\boldsymbol{x}_{\mathcal{P}_{n_k}, i}$ denote the $i$-th sample point in the subset $\mathcal{P}_{n_k}$. The point set minimizing the surrogate function $g(\mathcal{P}_n \mid \mathcal{P}_n^{(t)})$ is given by:  
\begin{equation}
 \begin{aligned}
\boldsymbol{x}_{\mathcal{P}_{n_k,i}} = &M_{n_k,i}(\mathcal{P}_{n}^{(t)}; \{\boldsymbol{y}_m\}_{m=1}^N)\\
  = &\left(\frac{1}{N}\sum_{m=1}^N\frac{1}{\| \boldsymbol{x}_{\mathcal{P}_{n_k,i}}^{(t)}-\boldsymbol{y}_m \|_2}\right)^{-1}\left(\frac{1}{N}\sum_{m=1}^N
 \frac{\boldsymbol{y}_m}{\|\boldsymbol{x}^{(t)}_{\mathcal{P}_{n_k,i}}-\boldsymbol{y}_m  \|_2}\right.\\
 &\left.+(1-\lambda)\frac{1}{n_k}\sum_{j\in\mathcal{I}_{n_k}\atop j\neq i}\frac{\boldsymbol{x}^{(t)}_{\mathcal{P}_{n_k,i}}-\boldsymbol{x}^{(t)}_j}{\|\boldsymbol{x}^{(t)}_{\mathcal{P}_{n_k,i}}-\boldsymbol{x}^{(t)}_j\|_2} +\lambda\frac{1}{n}\sum_{j=1\atop j\neq i}^n \frac{\boldsymbol{x}^{(t)}_{\mathcal{P}_{n_k,i}}-\boldsymbol{x}^{(t)}_j}{\|\boldsymbol{x}^{(t)}_{\mathcal{P}_{n_k,i}}-\boldsymbol{x}^{(t)}_j\|_2}\right),
\label{closedform1}
\end{aligned}
\end{equation}
where $ i =  1, \cdots, n_k,\ k=1,\dots,K.$
\end{theorem}  

According to Theorem \ref{Theorem6} and the framework of the MM algorithm, the construction algorithm for sliced space-filling design with mixtures based on the hybrid criterion  $E^\lambda(F, F_{\mathcal{P}_n})$ is presented in Algorithm \ref{algOne-shot}. The convergence of the iterative process of the MM algorithm is established in Theorem \ref{Theorem7}.  

\IncMargin{1em}
\begin{algorithm}[H]
\caption{One-shot algorithm for sliced space-filling design with mixtures}
\label{algOne-shot}
    \SetAlgoNoLine
 Generate uniform samples $ \{\boldsymbol{y}_m\}_{m=1}^N $ in the mixture experimental region $ \mathcal{X} $, where $ N \gg n $, and set $ t = 0 $;\\  
 Generate the initial point set $ \mathcal{P}_n^{(0)} = \{\boldsymbol{x}_i^{(0)} = \boldsymbol{y}_{u_i} + \boldsymbol{\varepsilon}_i\}_{i=1}^n $, where $ u_1, \ldots, u_n $ are sampled without replacement from $ \{1, \ldots, N\} $, and $ \boldsymbol{\varepsilon}_1, \ldots, \boldsymbol{\varepsilon}_n \in U[-\tau, \tau]^p $;\\  
 Randomly partition $ \mathcal{P}_n^{(0)} $ into $ K $ subsets $ \mathcal{P}_{n_1}^{(0)}, \ldots, \mathcal{P}_{n_K}^{(0)} $ with sample sizes $ n_1, \ldots, n_K $, respectively;\\  
\Repeat{$ \mathcal{P}_n^{(t)} $ converges}{
  In parallel, compute $ \boldsymbol{x}_{\mathcal{P}_{n_k}, i}^{(t+1)} = M_{n_k, i}\left(\mathcal{P}_n^{(t)}; \{\boldsymbol{y}_m\}_{m=1}^N\right) $ for $ i = 1, \ldots, n_k $ and $ k = 1, \ldots, K $;\\
 Update $ \mathcal{P}_{n_k}^{(t+1)} = \{\boldsymbol{x}_{\mathcal{P}_{n_k}, i}^{(t+1)}\}_{i=1}^{n_k} $ for $ k = 1, \ldots, K $, and set $ \mathcal{P}_n^{(t)} = \bigcup_{k=1}^K \mathcal{P}_{n_k}^{(t)} $, then $ t \gets t + 1 $;\\  } 
 \Return the converged point sets $ \mathcal{P}_{n_1}^*, \ldots, \mathcal{P}_{n_K}^* $, and $ \mathcal{P}_n^* = \bigcup_{k=1}^K \mathcal{P}_{n_k}^* $.  
\end{algorithm}
\DecMargin{1em}

\begin{theorem}\label{Theorem7}  
Suppose $ \mathcal{X} \subseteq T^p $ is a closed convex set. Given training samples $ \{\boldsymbol{y}_m\}_{m=1}^N \subseteq \mathcal{X} $, for any distinct initial point sets $ \mathcal{P}_n^{(0)} \subseteq \mathcal{X} $, the limit point set of the sequence $ \{\mathcal{P}_n^{(t)}\}_{t=0}^\infty $ in Algorithm $\ref{algOne-shot}$ converges to a stationary solution of the optimization problem $\eqref{ophy}$.  
\end{theorem} 

Considering the sequential form of method in section \ref{method3}, the optimization problem \eqref{seqmix} can be approximated as follows:
\begin{equation}\label{mixseqpractice}
\begin{aligned}
\underset{\mathcal{P}_{n_k} \subseteq \mathcal{X}}{ \textup{min}}\  h^s(\mathcal{P}_{n_k})=& \frac{2}{n_kN}\sum_{m=1}^N\sum_{i\in\mathcal{I}_{n_{k}}}\|\boldsymbol{x}_i-\boldsymbol{y}_m\|_2-
\frac{1}{n_k^2}\left(1-\lambda_k\frac{ n_c}{n_c+n_k}\right)\sum_{i,j\in\mathcal{I}_{n_{k}}}\|\boldsymbol{x}_i-\boldsymbol{x}_j\|_2\\
&-\frac{2\lambda_k}{(n_c+n_k)n_k}\sum_{i\in\mathcal{I}_{n_{k}}}\sum_{j\in\mathcal{I}_{n_{c}}}\|\boldsymbol{x}_i-\boldsymbol{x}_j\|_2,
\end{aligned}
\end{equation}

Objective function $h ^ s (\mathcal {P} _ {n}) $ in the current iteration point set $\mathcal {P} _ {n} ^ {(t)} $ agent function is given in Theorem \ref{Theorem8}, the proof ideas with Theorem \ref{Theorem6}.

\begin{theorem}\label{Theorem8}
Suppose $\{\boldsymbol{y}_m\}_{m=1}^N$ and $\mathcal{P}_{n_k}^{(t)} \subseteq \mathcal{X}$. Define the following function:  
\begin{equation*}
\begin{split}
&g^s(\mathcal{P}_{n_k} \mid \mathcal{P}_{n_k}^{(t)})\\ =\ & \frac{1}{n_k N} \sum_{m=1}^N \sum_{i \in \mathcal{I}_{n_k}} \left( \frac{\|\boldsymbol{x}_i - \boldsymbol{y}_m\|_2^2}{\|\boldsymbol{x}_i^{(t)} - \boldsymbol{y}_m\|_2} + \|\boldsymbol{x}_i^{(t)} - \boldsymbol{y}_m\|_2 \right) \\
&\quad - \frac{1}{n_k^2} \left( 1 - \lambda_k \frac{n_c}{n_c + n_k} \right) \sum_{\substack{i,j \in \mathcal{I}_{n_k} \\ i \neq j}} \left[ \|\boldsymbol{x}_i^{(t)} - \boldsymbol{x}_j^{(t)}\|_2 + \frac{2(\boldsymbol{x}_i^{(t)} - \boldsymbol{x}_j^{(t)})^T (\boldsymbol{x}_i - \boldsymbol{x}_i^{(t)})}{\|\boldsymbol{x}_i^{(t)} - \boldsymbol{x}_j^{(t)}\|_2} \right] \\
&\quad - \frac{2\lambda_k}{(n_c + n_k) n_k} \sum_{i \in \mathcal{I}_{n_k}} \sum_{j \in \mathcal{I}_{n_c}} \left[ \|\boldsymbol{x}_i^{(t)} - \boldsymbol{x}_j^{(t)}\|_2 + \frac{(\boldsymbol{x}_i^{(t)} - \boldsymbol{x}_j^{(t)})^T (\boldsymbol{x}_i - \boldsymbol{x}_i^{(t)})}{\|\boldsymbol{x}_i^{(t)} - \boldsymbol{x}_j^{(t)}\|_2} \right],
\end{split}
\end{equation*}  
Then, $g^s(\mathcal{P}_{n_k} \mid \mathcal{P}_{n_k}^{(t)})$ is a surrogate function of the objective function $h^s(\mathcal{P}_{n_k})$ in \eqref{mixseqpractice} at the current iteration $\mathcal{P}_{n_k}^{(t)}$, satisfying:  
$(i)$ $g^s(\mathcal{P}_{n_k} \mid \mathcal{P}_{n_k}^{(t)}) \geq h^s(\mathcal{P}_{n_k})$ for all $\mathcal{P}_{n_k} \subseteq \mathcal{X}$;  
$(ii)$ $g^s(\mathcal{P}_{n_k}^{(t)} \mid \mathcal{P}_{n_k}^{(t)}) = h^s(\mathcal{P}_{n_k}^{(t)})$.  
Further, let $\boldsymbol{x}_{\mathcal{P}_{n_k}, i}$ denote the $i$-th sample point in the subset $\mathcal{P}_{n_k}$. The point set minimizing $g^s(\mathcal{P}_{n_k} \mid \mathcal{P}_{n_k}^{(t)})$ is given by the following for $i = 1, \ldots, n_k$:  
\begin{equation}\label{closedform2}
\begin{split}
\boldsymbol{x}_{\mathcal{P}_{n_k}, i} &= M_{n_k, i}^s\left( \mathcal{P}_n^{(t)}; \{\boldsymbol{y}_m\}_{m=1}^N \right) \\
&= \left( \frac{1}{N} \sum_{m=1}^N \frac{1}{\|\boldsymbol{y}_m - \boldsymbol{x}_{\mathcal{P}_{n_k}, i}^{(t)}\|_2} \right)^{-1} \left( \frac{1}{N} \sum_{m=1}^N \frac{\boldsymbol{y}_m}{\|\boldsymbol{x}_{\mathcal{P}_{n_k}, i}^{(t)} - \boldsymbol{y}_m\|_2} \right. \\
&\quad \left. + \left( 1 - \lambda_k \frac{n_c}{n_c + n_k} \right) \frac{1}{n_k} \sum_{\substack{j \in \mathcal{I}_{n_k} \\ j \neq i}} \frac{\boldsymbol{x}_{\mathcal{P}_{n_k}, i}^{(t)} - \boldsymbol{x}_j^{(t)}}{\|\boldsymbol{x}_{\mathcal{P}_{n_k}, i}^{(t)} - \boldsymbol{x}_j^{(t)}\|_2} +\right. \\ &\left. \left( \lambda_k \frac{n_c}{n_c + n_k} \right) \frac{1}{n_c} \sum_{j \in \mathcal{I}_{n_c}} \frac{\boldsymbol{x}_{\mathcal{P}_{n_k}, i}^{(t)} - \boldsymbol{x}_j^{(t)}}{\|\boldsymbol{x}_{\mathcal{P}_{n_k}, i}^{(t)} - \boldsymbol{x}_j^{(t)}\|_2} \right).
\end{split}
\end{equation}  
\end{theorem}

According to Theorem \ref{Theorem8} and the basic framework of the Majorization-Minimization (MM) algorithm, the sequential construction algorithm for sliced space-filling design with mixtures based on the hybrid criterion $ E^{\lambda_k}(F, F_{\mathcal{P}_n}) $ is presented in Algorithm \ref{algSequential}. The convergence of the MM algorithm is established in Theorem \ref{Theorem9}.  

\begin{theorem}\label{Theorem9}
Suppose $\mathcal{X} \subseteq T^p$ is a closed convex set. Given training samples $\{\boldsymbol{y}_m\}_{m=1}^N \subseteq \mathcal{X}$, for any distinct initial point sets $\mathcal{P}_{n_k}^{(0)} \subseteq \mathcal{X}$, the limit point set of the sequence $\{\mathcal{P}_{n_k}^{(t)}\}_{t=0}^\infty$ in Algorithm $2$ converges to a stationary solution of the optimization problem \eqref{mixseqpractice}.  
\end{theorem}  
\IncMargin{1em}
\begin{algorithm}[H]
\caption{Sequential algorithm for sliced space-filling design with mixtures}
\label{algSequential}
    \SetAlgoNoLine
 Generate uniform samples $ \{\boldsymbol{y}_m\}_{m=1}^N $ in the mixture experimental region $ \mathcal{X} $, where $ N \gg n $;  \\
 Initialize $ n_c = t = 0 $ and $ \mathcal{P}_{n_c} = \emptyset $;  \\
\For{$ k = 1, \ldots, K $} {
   Generate the initial point set $ \mathcal{P}_{n_k}^{(0)} = \{\boldsymbol{x}_i^{(0)} = \boldsymbol{y}_{u_i} + \boldsymbol{\varepsilon}_i\}_{i=1}^{n_k} \subseteq \mathcal{X} $, where $ u_1, \ldots, u_{n_k} $ are sampled without replacement from $ \{1, \ldots, N\} $, and $ \boldsymbol{\varepsilon}_1, \ldots, \boldsymbol{\varepsilon}_{n_k} \in U[-\tau, \tau]^p $;  \\
    \Repeat{$ \mathcal{P}_{n_k}^{(t)} $ converges, denoting the converged set as $ \mathcal{P}_{n_k}^* $}{    
         In parallel, compute $ \boldsymbol{x}_{\mathcal{P}_{n_k}, i}^{(t+1)} = M_{n_k, i}^s\left(\mathcal{P}_{n_k}^{(t)}; \{\boldsymbol{y}_m\}_{m=1}^N\right) $ for $ i = 1, \ldots, n_k $;  \\
         Update $ \mathcal{P}_{n_k}^{(t+1)} = \{\boldsymbol{x}_{\mathcal{P}_{n_k}, i}^{(t+1)}\}_{i=1}^{n_k} $ and set $ t \gets t + 1 $;  \\
    }
    Update $ \mathcal{P}_{n_c} \gets \mathcal{P}_{n_c} \cup \mathcal{P}_{n_k}^* $ and $ n_c \gets n_c + n_k $;\\
    } 
  Return the converged point sets $ \mathcal{P}_{n_1}^*, \ldots, \mathcal{P}_{n_K}^* $, and $ \mathcal{P}_n^* = \bigcup_{k=1}^K \mathcal{P}_{n_k}^* $.  
\end{algorithm}
\DecMargin{1em}

\begin{remark} 
The tuning parameters $ \lambda $ and $ \lambda_k $ aim to balance the importance between the representativeness of the full point set and that of the subsets, and their selection depends on specific application scenarios. In the absence of prior information, based on the iteratives $\eqref{closedform1}$ and $\eqref{closedform2}$ in Theorems $\ref{Theorem6}$ and $\ref{Theorem8}$ it is recommended to set $ \lambda = 0.5 $ in Algorithm $\ref{algOne-shot}$. For the sequential sampling in Algorithm $\ref{algSequential}$, to ensure both the full point set and subsets, $ \lambda_k $ can be set as: $\lambda_k = \frac{n_c + n_k}{2n_c} \quad (k \geq 2)$
$( $When $ k = 1 $, $ n_c = 0 $ and $ \mathcal{I}_{n_c} = \emptyset $,  the value of $ \lambda_1 $ has no impact on the optimization problem $\eqref{mixseqpractice}$.$)$ The underlying idea is that during the iterative adjustment process of each point, equal weights are assigned to the average adjustment directions corresponding to the two sets indexed by $ j $. All examples in this paper adopt this scheme. If the representativeness of the full point set is more important, $ \lambda, \lambda_k \in [0.7, 1] $ can be chosen; if the representativeness of the subsets is more important, $ \lambda, \lambda_k \in [0, 0.3] $ can be chosen.  
\end{remark}  

\begin{remark}
In the first step of Algorithms $\ref{algOne-shot}$ and $\ref{algSequential}$, the training samples $ \{\boldsymbol{y}_m\}_{m=1}^N $ can be generated using sampling methods such as the inverse transform method $\citep{Wang1990}$, Gibbs sampling $\citep{xiong2024}$, and hit-and-run sampling $\citep{hitandrunr}$. $ \tau $ is a sufficiently small positive number.
\end{remark}  

\begin{remark}  
When computing $ \boldsymbol{x}_{\mathcal{P}_{n_k}, i}^{(t+1)} $ in each iteration, a subset of $ N_s $ samples can be randomly drawn from $ \{\boldsymbol{y}_m\}_{m=1}^N $ for the calculation to improve computational efficiency.  
\end{remark}  

\section{ Numerical Experiments}

This section presents several numerical examples to compare the performance of mixed fragment designs constructed by different algorithms in terms of uniformity and modeling prediction, including the inverse transformation method (\textbf{ITM}), coset decomposition method (\textbf{CDM}), and fast flexible space-filling design (\textbf{FFFD}) that are already available in the literature mentioned in section \ref{sec1}. And the algorithms newly proposed in section \ref{sec3} of this paper based on the energy distance criterion, including sequentially generate subdesign without balance parameter (\textbf{SeqM}), partition the full point set (\textbf{ParM}) and one-shot and sequential algorithms based on hybrid criterion (denoted as \textbf{MHED} and \textbf{SeqHED} respectively). To ensure that most methods are applicable, the following mainly considers the experimental region of the mixture components as the standard simplex. 

\subsection{ Uniformity Comparison of sliced space-filling design with mixtures}  
\begin{example}\label{ex1}
    Generate sliced space-filling design with mixtures in the following regions$:$  
\begin{enumerate}[$(i)$]
  \item  $\mathcal{X}^1_\mathcal{Z}=\{(\boldsymbol{x}^T, z)^T\ |\ \boldsymbol{x}\in T_3,\ z\in\{1,2,3\}\}$;
  \item $\mathcal{X}^2_\mathcal{Z}=\{(\boldsymbol{x}^T, \boldsymbol{z})^T\ |\ \boldsymbol{x}=(x_1,x_2,x_3)^T,0.1\leq x_1\leq 0.8,\ 0.05\leq x_2\leq 0.6,\ 0.15\leq x_3\leq 0.7,\sum_{i=1}^3x_i=1,\ \boldsymbol{z}=(z_1,z_2)^T\in\{-1,1\}\}$.
\end{enumerate}
\end{example}

Generate uniform random numbers $ \{\boldsymbol{u}_t\}_{t=1}^N $ in the experimental region. For the full design $ \mathcal{P}_n $ and subdesigns $ \mathcal{P}_{n_k} $ produced by each method, we compute the following uniformity evaluation criteria:  
 Estimation errors of the mean ($ \Delta_\mu $) and standard deviation ($ \Delta_\sigma $) relative to the uniform distribution in the target region; Root sean squared Distance (RMSD); Maximum packing distance (MaD); Minimum separation distance (MiD).  
\begin{equation*}
\begin{aligned}
& \Delta_{\mu}(\mathcal{P}_n)=\|\boldsymbol{\hat{\mu}}(\mathcal{P}_n)-\boldsymbol{\mu}\|_2+\sum_{k=1}^K\|\boldsymbol{\hat{\mu}}(\mathcal{P}_{n_k})-\boldsymbol{\mu}\|_2,\\
& \Delta_{\sigma}(\mathcal{P}_n)=\|\boldsymbol{\hat{\sigma}}(\mathcal{P}_n)-\boldsymbol{\sigma}\|_2+\sum_{k=1}^K\|\boldsymbol{\hat{\sigma}}(\mathcal{P}_{n_k})-\boldsymbol{\sigma}\|_2,\\ &\mathrm{RMSD}(\mathcal{P}_n)=\sqrt{\frac{1}{N}\sum_{t=1}^N\underset{\boldsymbol{x}\in\mathcal{P}_n}{\min}\|\boldsymbol{u}_t
 -\boldsymbol{x}\|_2^2}+\sum_{k=1}^K\sqrt{\frac{1}{N}\sum_{t=1}^N\underset{\boldsymbol{x}\in\mathcal{P}_{n_k}}{\min}\|\boldsymbol{u}_t
 -\boldsymbol{x}\|_2^2},\\
& \mathrm{MaD}(\mathcal{P}_n)=\max_{t\in\{1,2,\dots,N\}}\underset{\boldsymbol{x}\in\mathcal{P}_n}{\min}\|\boldsymbol{u}_t-\boldsymbol{x}\|_2
+\sum_{k=1}^K\max_{t\in\{1,2,\dots,N\}}\underset{\boldsymbol{x}\in\mathcal{P}_{n_k}}{\min}\|\boldsymbol{u}_t-\boldsymbol{x}\|_2,\\
& \mathrm{MiD}(\mathcal{P}_n)=\underset{\boldsymbol{x}_i,\boldsymbol{x}_j\in\mathcal{P}_n\atop \boldsymbol{x}_i\neq\boldsymbol{x}_j}{\min}\|\boldsymbol{x}_i-\boldsymbol{x}_j\|_2
+\sum_{k=1}^K\underset{\boldsymbol{x}_i,\boldsymbol{x}_j\in\mathcal{P}_{n_k}\atop \boldsymbol{x}_i\neq\boldsymbol{x}_j}{\min}\|\boldsymbol{x}_i-\boldsymbol{x}_j\|_2.
 \end{aligned}
\end{equation*}

For the target region $\mathcal{X}_\mathcal{Z}^1$, the process variable $z$ has 3 levels, so we consider the scenario where the full mixture component design is partitioned into 3 slices. All methods are applicable in this scenario.  
For the ITM method, we first generate a maximin Latin hypercube design using the R package \texttt{SLHD}, and then map it to $\mathcal{X}_\mathcal{Z}^1$ via a transformation function. The FFFM method generates sliced space-filling design with mixtures using the space-filling design module in JMP 18. Details of design generation for other methods are described in this paper and related literature, with parameters set to their default values.  
The process of generating sliced space-filling design with mixtures via each method is repeated 50 times, and the uniformity criterion values for these designs are computed; the results are shown in Figure\ref{fig2}. The scatter plots of these designs are shown in the supplementary material.

\begin{figure}[htb]
    \centering
    \begin{subfigure}[b]{0.32\textwidth}
        \includegraphics[width=\textwidth]{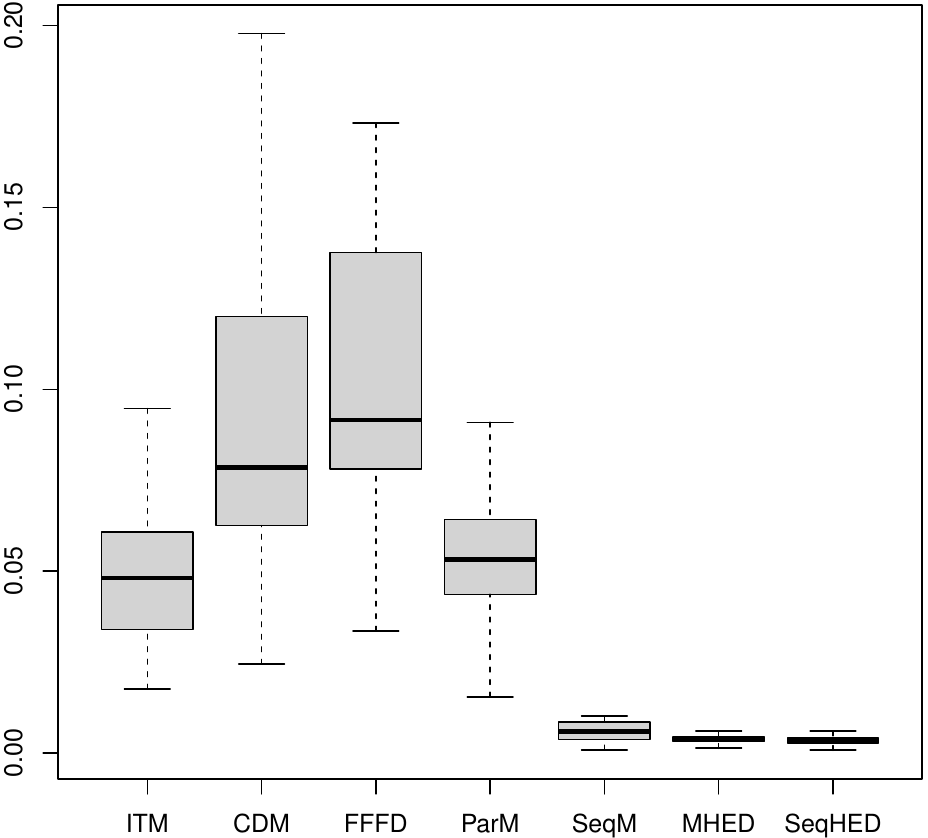}
        \caption{ $\Delta_{\mu}$}
     \end{subfigure}
    \begin{subfigure}[b]{0.32\textwidth}
        \includegraphics[width=\textwidth]{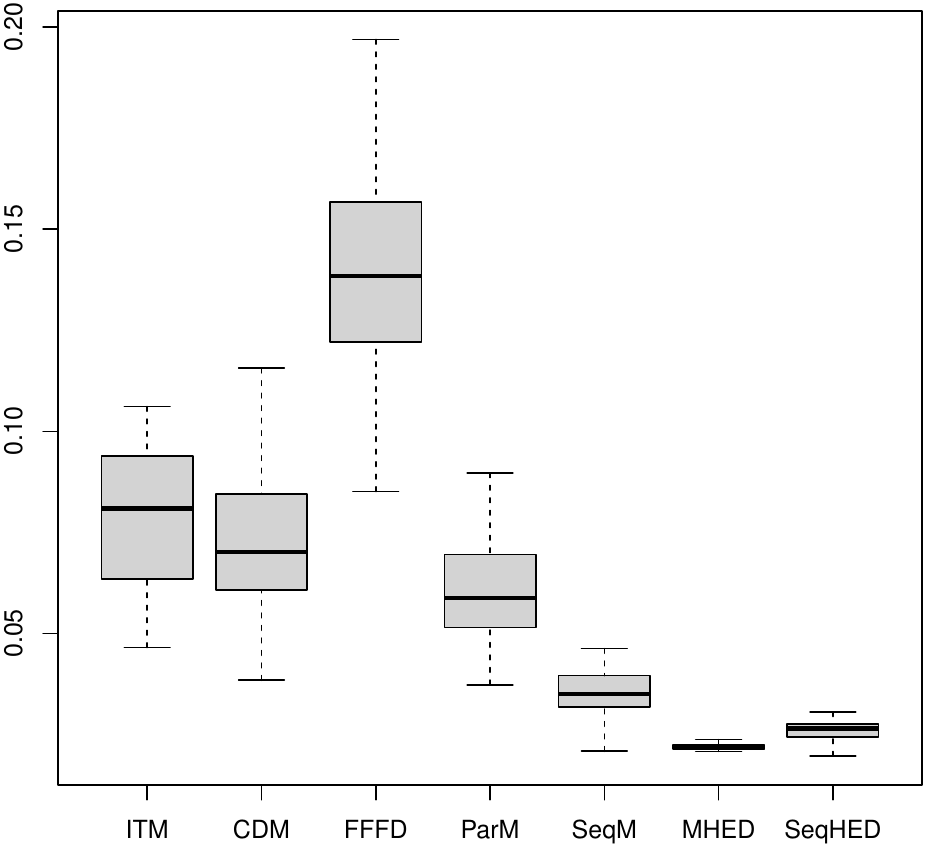}
        \caption{ $\Delta_{\sigma}$}
  \end{subfigure}
  \begin{subfigure}[b]{0.32\textwidth}
        \includegraphics[width=\textwidth]{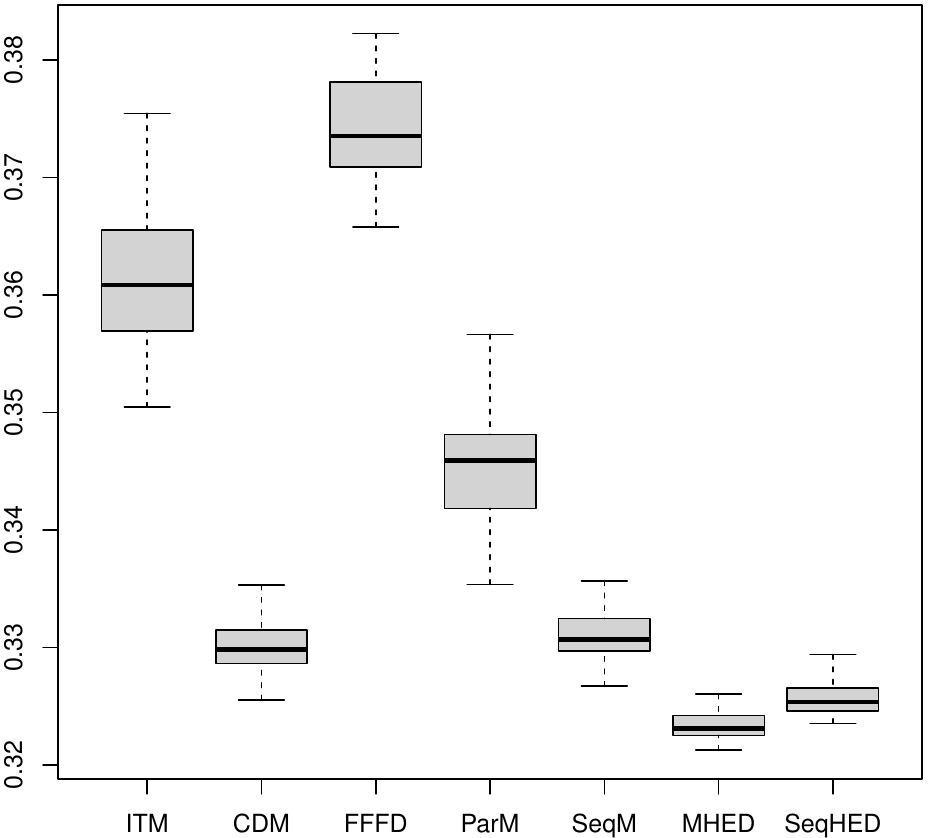}
        \caption{RMSD}
  \end{subfigure}
    \begin{subfigure}[b]{0.32\textwidth}
        \includegraphics[width=\textwidth]{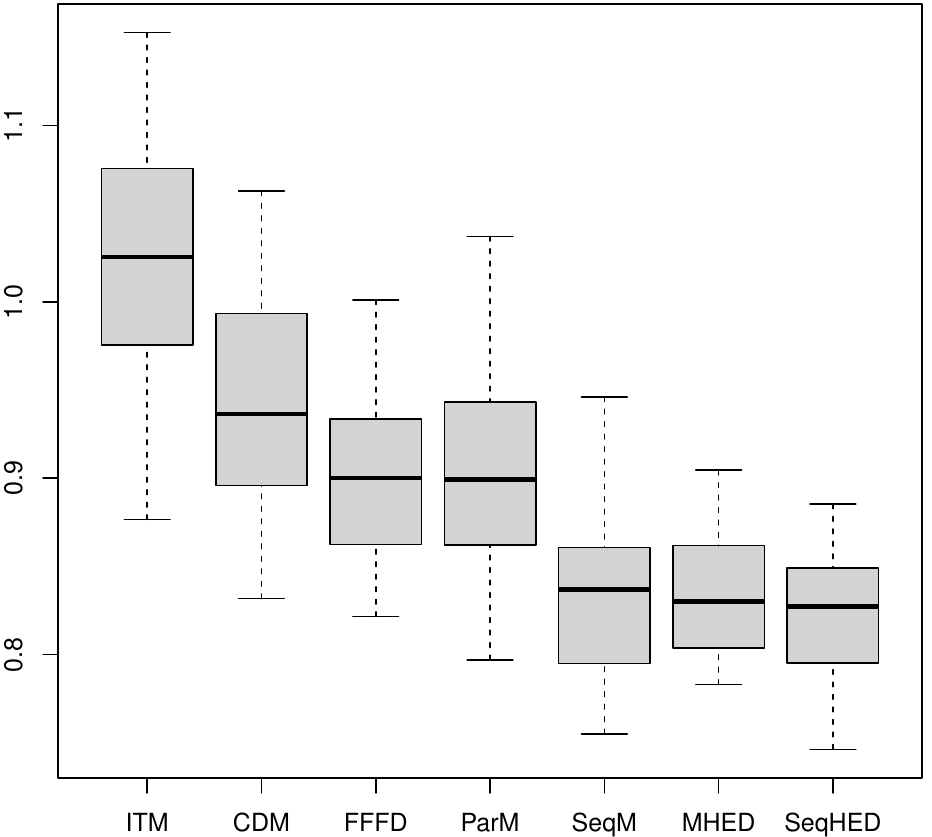}
        \caption{MaD}
  \end{subfigure}
    \begin{subfigure}[b]{0.32\textwidth}
        \includegraphics[width=\textwidth]{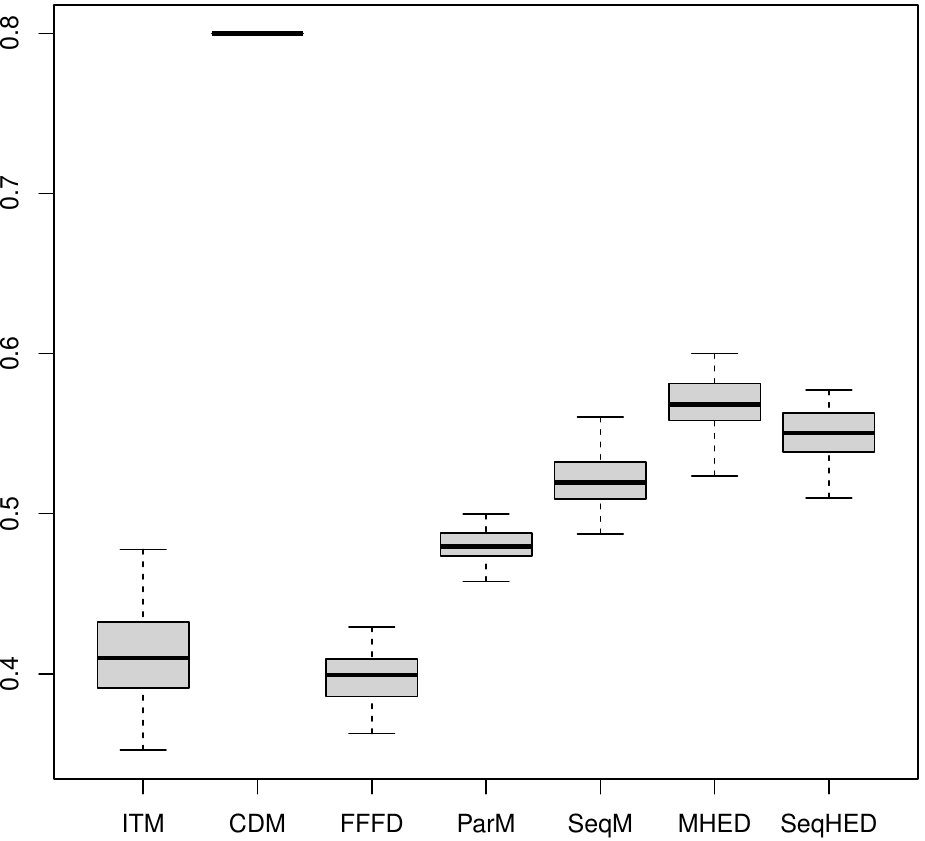}
        \caption{MiD}
  \end{subfigure}
 \caption{\ Comparison of uniformity for sliced space-filling design with mixtures generated by different methods on $\mathcal{X}^1_\mathcal{Z}$, Where larger MiD criterion values are preferred and smaller values are better for other criteria.}
\label{fig2}
\end{figure}

\begin{figure}[htb]
    \centering
    \begin{subfigure}[b]{0.32\textwidth}
        \includegraphics[width=\textwidth]{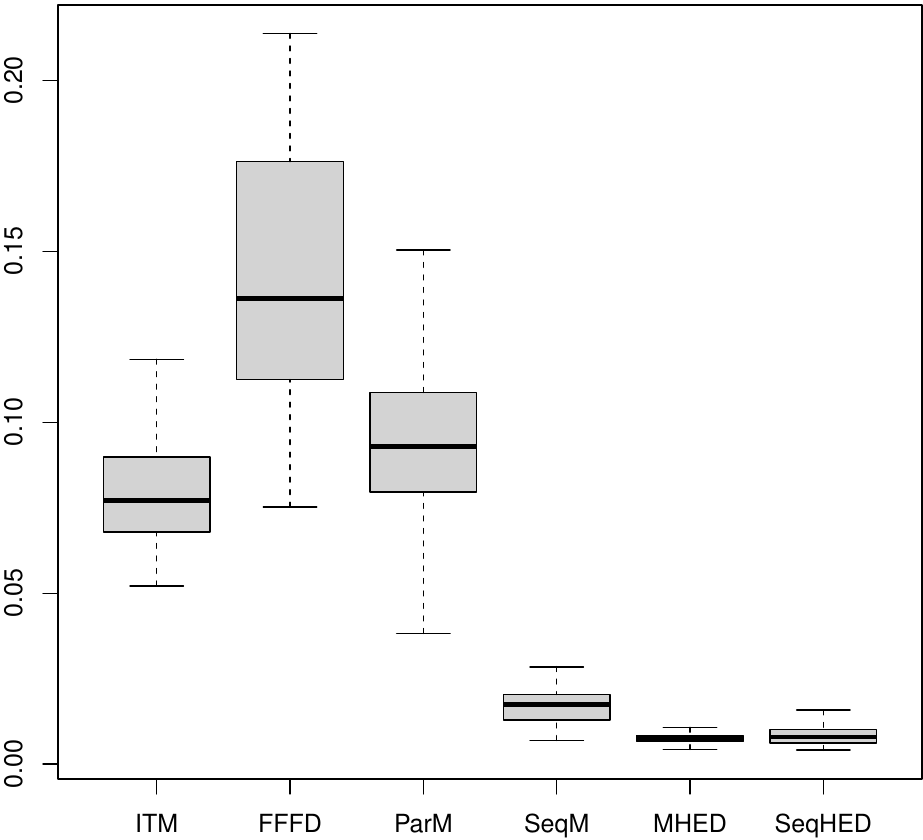}
        \caption{ $\Delta_{\mu}$}
     \end{subfigure}
    \begin{subfigure}[b]{0.32\textwidth}
        \includegraphics[width=\textwidth]{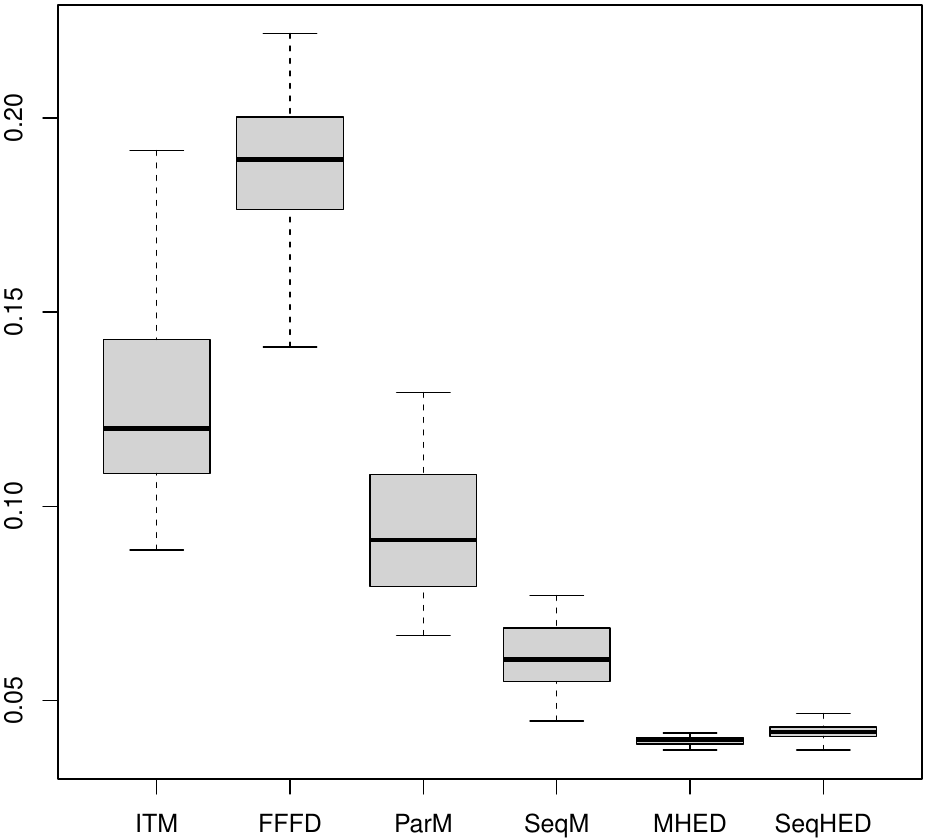}
        \caption{ $\Delta_{\sigma}$}
  \end{subfigure}
  \begin{subfigure}[b]{0.32\textwidth}
        \includegraphics[width=\textwidth]{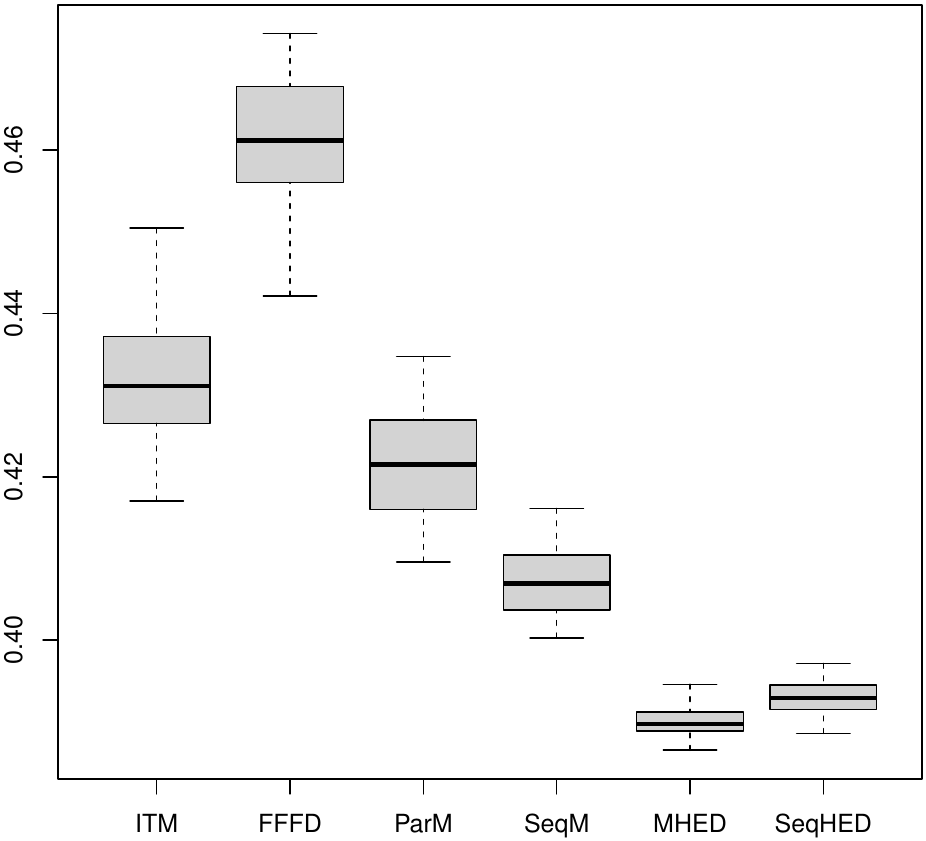}
        \caption{RMSD}
  \end{subfigure}
    \begin{subfigure}[b]{0.32\textwidth}
        \includegraphics[width=\textwidth]{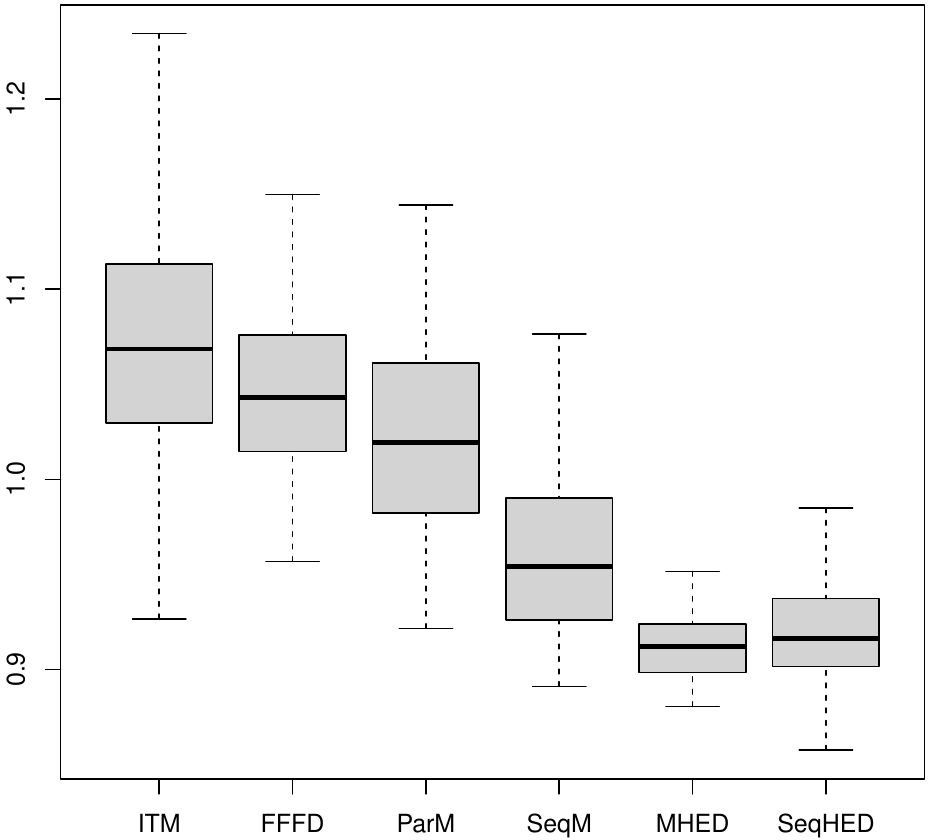}
        \caption{MaD}
  \end{subfigure}
    \begin{subfigure}[b]{0.32\textwidth}
        \includegraphics[width=\textwidth]{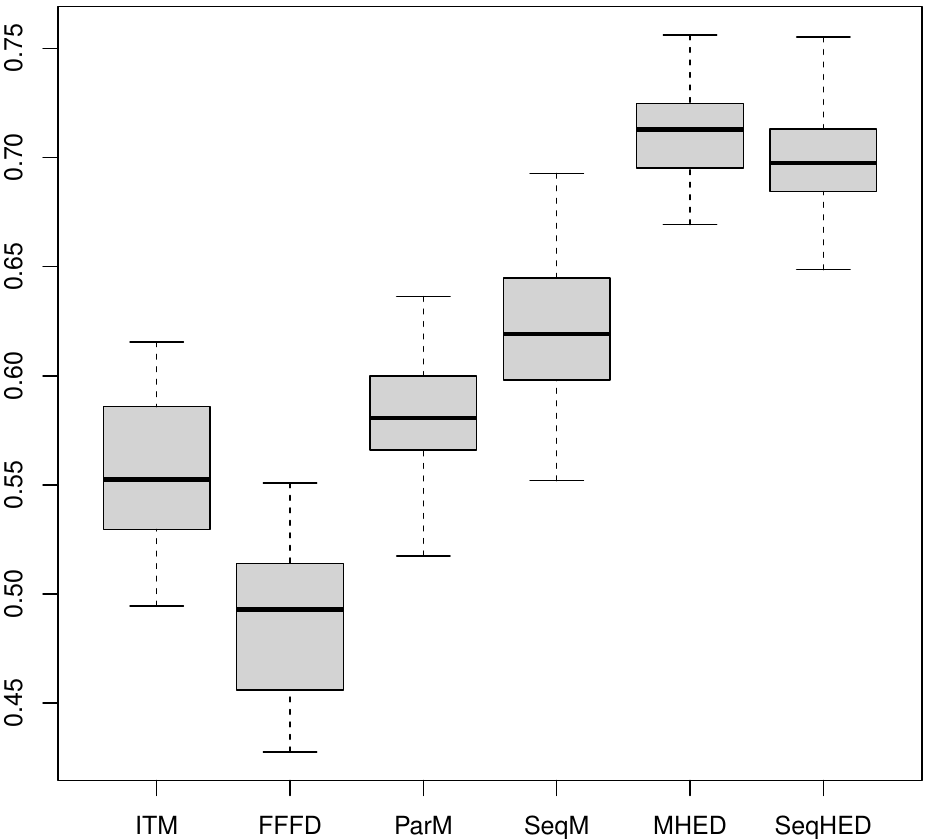}
        \caption{MiD}
  \end{subfigure}
 \caption{\ Comparison of uniformity for sliced space-filling design with mixtures generated by different methods on $\mathcal{X}^2_\mathcal{Z}$, where larger MiD criterion values are preferred and smaller values are better for other criteria.}
\label{fig3}
\end{figure}

As observed in Figure\ref{fig2}(e), the CDM method outperforms MHED under the MiD criterion. The main reason is that CDM generates designs based on a lattice point structure, so the design points in both the full design and subdesigns exhibit excellent pairwise separation properties.  
Besides this, the proposed MHED method performs best, followed by SeqHED, which is far superior in uniformity to other methods. However, from the perspective of sequential experimental needs in practice, SeqHED is more flexible. The sliced space-filling design with mixtures generated by ITM (by transforming well-uniformed partition designs on the hypercube) are not necessarily uniform; in particular, the uniformity of subdesigns may be poor.  
Although the ParM method ensures optimal uniformity for the full design, the uniformity of subdesigns generated by partitioning may be poor when the sample size is small. Additionally, solving the optimization problem in \eqref{par2} is relatively challenging. Thus, compared to other methods based on the energy distance criterion proposed in this paper, the performance of ParM is slightly inferior. 

For the target region $\mathcal{X}_\mathcal{Z}^2$, the process variable $\boldsymbol{z} = (z_1, z_2)^T$ has 4 level combinations, so we consider the scenario where the full design is partitioned into 4 slices. Since the CDM method is only applicable to the standard simplex case, we exclude this method from the comparison.  

The process of generating sliced space-filling design with mixtures via the remaining methods is repeated 50 times, and the relevant results are presented in Figure\ref{fig3}. The conclusions drawn from Figure\ref{fig3} are consistent with those from Figure\ref{fig2}: the uniformity of sliced space-filling design with mixtures constructed based on the energy distance criterion is significantly superior to that of existing methods.
\subsection{Prediction Performance}

In this subsection, we compare the prediction performance of the classical mixture design with process variables (\textbf{CMP}) introduced in section \ref{sec2.1} against the various sliced space-filling design with mixtures mentioned above. Additionally, we also consider two methods:
\textbf{randParM}: Generates a uniformly distributed full design and then randomly partitions it into subdesigns;
\textbf{comM}: Independently generates well-uniformed subdesigns and then aggregates them into a full design.

In some mixture experiments, to accelerate product development, experimenters often distribute experimental points of a mixture design across $ K $ laboratories (or process variables such as instruments, simulation platforms, etc.) for simultaneous experiments. After completing the experiments, data is either analyzed independently or aggregated for joint analysis. In such cases, different levels of the process variable correspond to the same response surface. We use a Gaussian process model to fit the experimental data from both the full design $ \mathcal{P}_n $ and its sub-designs $ \mathcal{P}_{n_k} $, and evaluate performance by computing the total root mean squared prediction error (RMSE) in the test set $ \{\boldsymbol{x}_m^*\}_{m=1}^M $:  
\begin{equation}\label{RMSEcri}
\text{RMSE}(\mathcal{P}_n) = \sqrt{\frac{1}{M} \sum_{m=1}^M \left( \hat{Y}_{\mathcal{P}_n}(\boldsymbol{x}_m^*) - Y(\boldsymbol{x}_m^*) \right)^2} + \sum_{k=1}^K \sqrt{\frac{1}{M} \sum_{m=1}^M \left( \hat{Y}_{\mathcal{P}_{n_k}}(\boldsymbol{x}_m^*) - Y(\boldsymbol{x}_m^*) \right)^2}
\end{equation}  
to assess the quality of the designs.  

\begin{example}\label{ex2}
Suppose that the real models are the three commonly used models in the mixture experiment: the quadratic model, the additional reciprocal model, and the Becker homogeneous model $\citep{Cornell2002,zhangli.2021}$. The experimental regions are all $\mathcal{X}^1_\mathcal{Z}$ in Example $\ref{ex1}$,\ $\forall\ z=1,2,3$: 
\begin{equation*}
\begin{aligned}
 &(i)\ f_1(x_1, x_2, x_3, z) =1-0.5x_1-0.5x_2+x_1^2+x_2^2+x_1x_2;
\\
& (ii)\ f_2(x_1, x_2, x_3, z)=1-0.5x_1-0.5x_2+\frac{1}{x_1+0.1}+\frac{1}{x_2+0.1};\\
&(iii)\   f_3(x_1, x_2, x_3, z) =0.5x_1+0.5x_2-0.5x_3+\min(x_1,x_2)+\min(x_1,x_3)+\min(x_2,x_3).
\end{aligned}
\end{equation*}
\end{example}

\begin{figure}[htb]
    \centering
    \begin{subfigure}[b]{0.31\textwidth}
        \includegraphics[width=\textwidth]{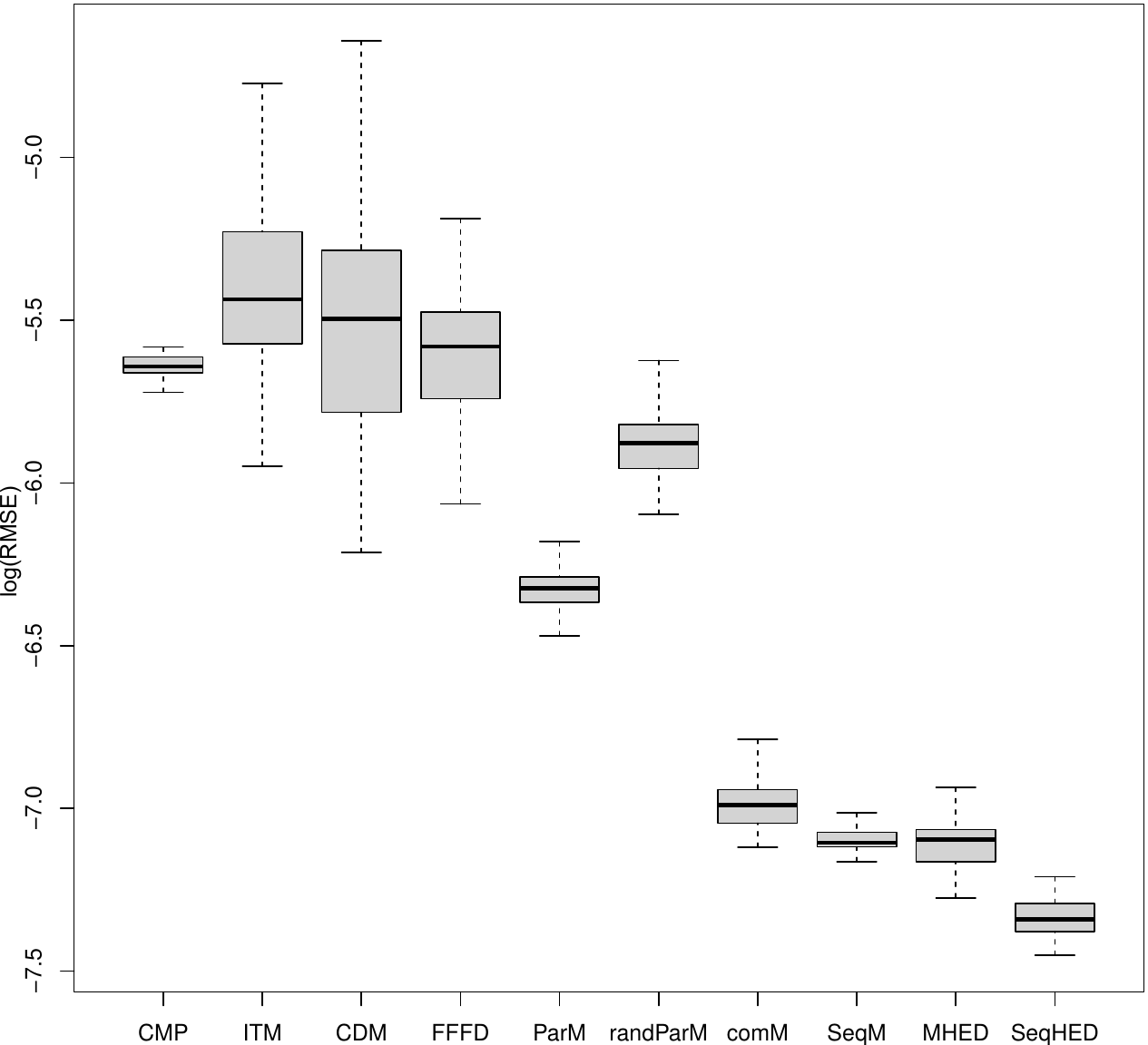}
        \caption{ $f_1$}
     \end{subfigure}
    \begin{subfigure}[b]{0.31\textwidth}
        \includegraphics[width=\textwidth]{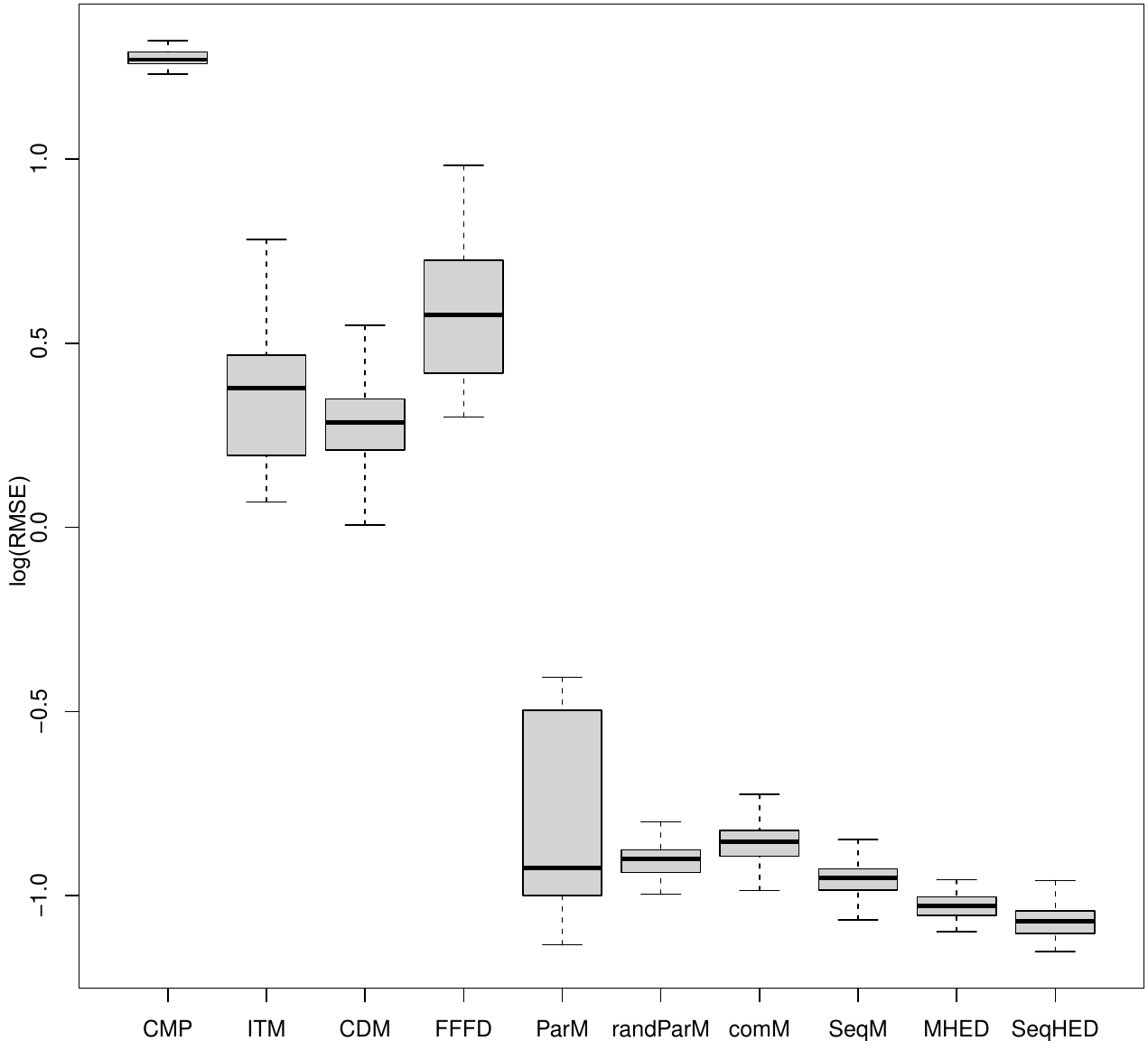}
        \caption{ $f_2$}
  \end{subfigure}
  \begin{subfigure}[b]{0.31\textwidth}
        \includegraphics[width=\textwidth]{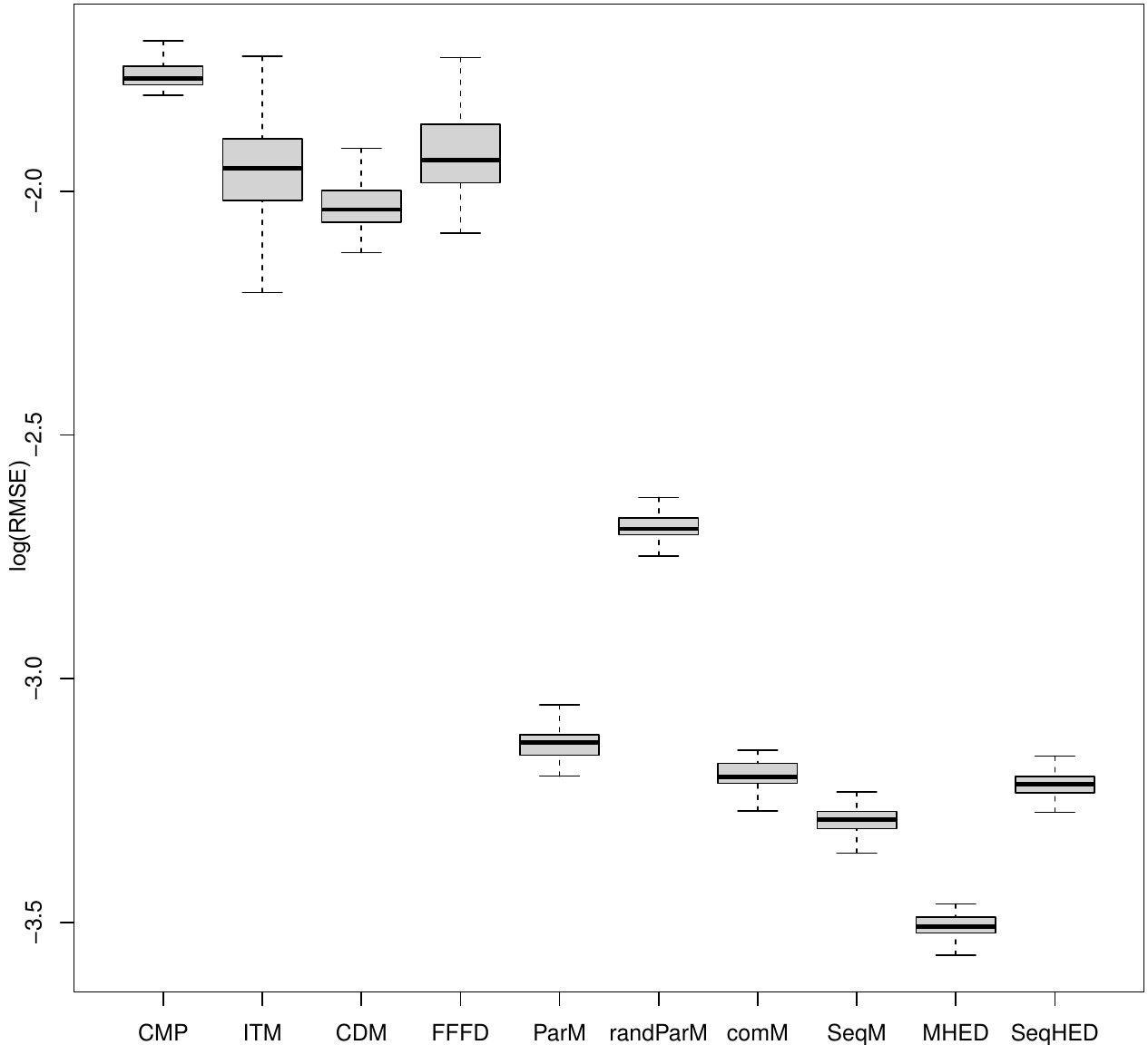}
        \caption{$f_3$}
  \end{subfigure}
 \caption{\ Comparison of modeling prediction errors based on different mixture design methods.}
\label{fig4}
\end{figure}

Set the number of experimental runs as $ n_1 = n_2 = n_3 = 15 $. For the CMP method, the mixture component can be specified as the $\{3,4\}$ simplex-lattice design \citep{Cornell2002}. Randomly generate $ M = 1000 $ uniformly distributed sample points in $ T_3 $ to serve as the test set. For the aforementioned models and design methods, each scenario is repeated 50 times, and the results are presented in Figure\ref{fig4}.  

From the results in Figure\ref{fig4}, it can be observed that the performance of all methods generally corresponds to the uniformity of the designs generated in Example \ref{ex1}. The six subsequent methods based on the energy distance criterion all significantly outperform the existing methods under the RMSE criterion. Among these, the MHED and SeqHED methods, which are based on the hybrid criterion, perform the best, and the SeqM method is also relatively robust. In contrast, the ParM method (which partitions a full design) and randParM exhibit somewhat inferior performance. The primary reason for this is the difficulty in solving the optimization problem \eqref{par2}. The solution method described in \cite{twinning}, which aims to enhance the efficiency of partitioning large-scale data, reduces the accuracy of the problem to a certain extent. Thus, more effective optimization algorithms warrant further investigation.  

\subsection{Application in drug combination experiment}
In many mixture problems, although the response surfaces of mixture components differ across different level combinations of process variables, they usually exhibit certain similarities, which may help improve prediction accuracy on each surface. We now consider the application of different sliced space-filling design with mixtures construction methods in drug mixture experiments.  

\begin{example} \label{ex3} 
$\cite{Ningdrug}$ studied the effects of the proportions $ \theta_1, \theta_2, \theta_3 $ of three drugs (AG490, U0126, and I-3-M) and the total dose $ C $ on the survival rate $ Y $ of normal cells and lung cancer cells, and established the following nonlinear regression model$:$  
\begin{equation}\label{drugmodel}
Y = \frac{1}{1 + \left( \frac{C}{IC_{50}(\boldsymbol{\theta})} \right)^{\gamma(\boldsymbol{\theta})}} + \varepsilon,
\end{equation}  
where $ \boldsymbol{\theta} = (\theta_1, \theta_2, \theta_3)^\mathrm{T} \in T_3 $, $ \varepsilon \sim N(0, \sigma^2) $, and $ IC_{50}(\boldsymbol{\theta}) $ and $ \gamma(\boldsymbol{\theta}) $ are second-order polynomial functions of $ \theta_1 $ and $ \theta_2 $. To compare the performance of slcied mixture designs in modeling and prediction, we take the fitted model for normal cells from $\cite{Ningdrug}$ as the true model, i.e., set $ \sigma = 10^{-4} $,  
\begin{equation*}
    IC_{50}(\boldsymbol{\theta}) = 117.11 - 15.71\theta_1 - 86.15\theta_2 + 79.55\theta_1^2 + 42.67\theta_2^2 - 33.75\theta_1\theta_2,
\end{equation*}
\begin{equation*}
    \gamma(\boldsymbol{\theta}) = 1.7 - 1.02\theta_1 + 0.41\theta_2 + 0.31\theta_1^2 - 0.74\theta_2^2 + 1.21\theta_1\theta_2.
\end{equation*}
The total dose of the drug mixture is $ C \in \{100, 200, 300, 400, 500, 600, 700\} $. Thus, we generate mixture designs in the experimental region:  
\begin{equation}
\mathcal{X}_\mathcal{C}^{\text{drug}} = \left\{ \boldsymbol{(\theta}^T, C)^T  \mid \boldsymbol{\theta} \in T_3,\ C \in \{100, 200, 300, 400, 500, 600, 700\} \right\}
\end{equation}  
and generate response values computed via model $\eqref{drugmodel}$.  
\end{example}

\begin{figure}[htb]
    \centering
    \begin{subfigure}[b]{0.32\textwidth}
        \includegraphics[width=\textwidth]{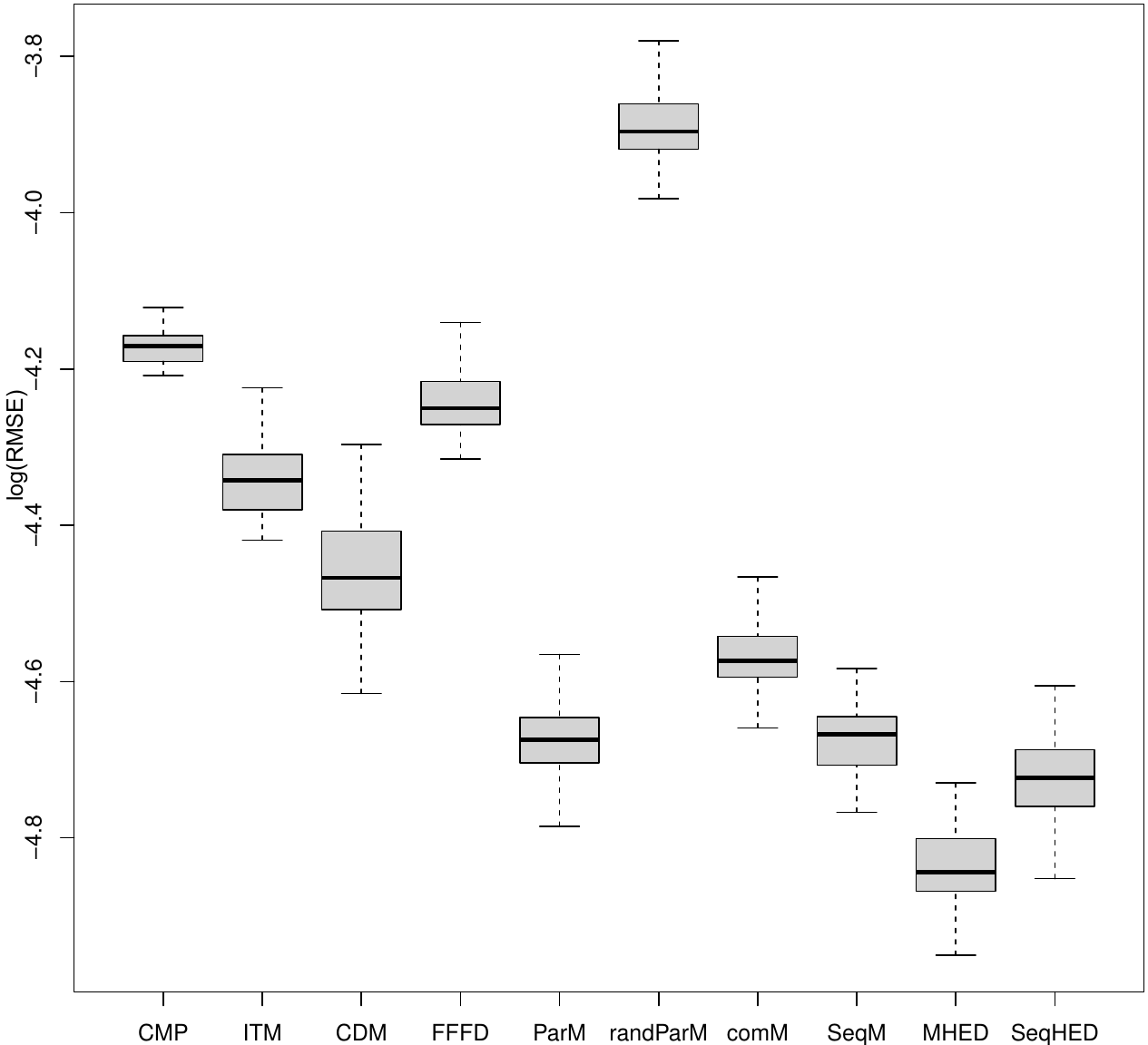}
        \caption{OGP}
     \end{subfigure}
    \begin{subfigure}[b]{0.32\textwidth}
        \includegraphics[width=\textwidth]{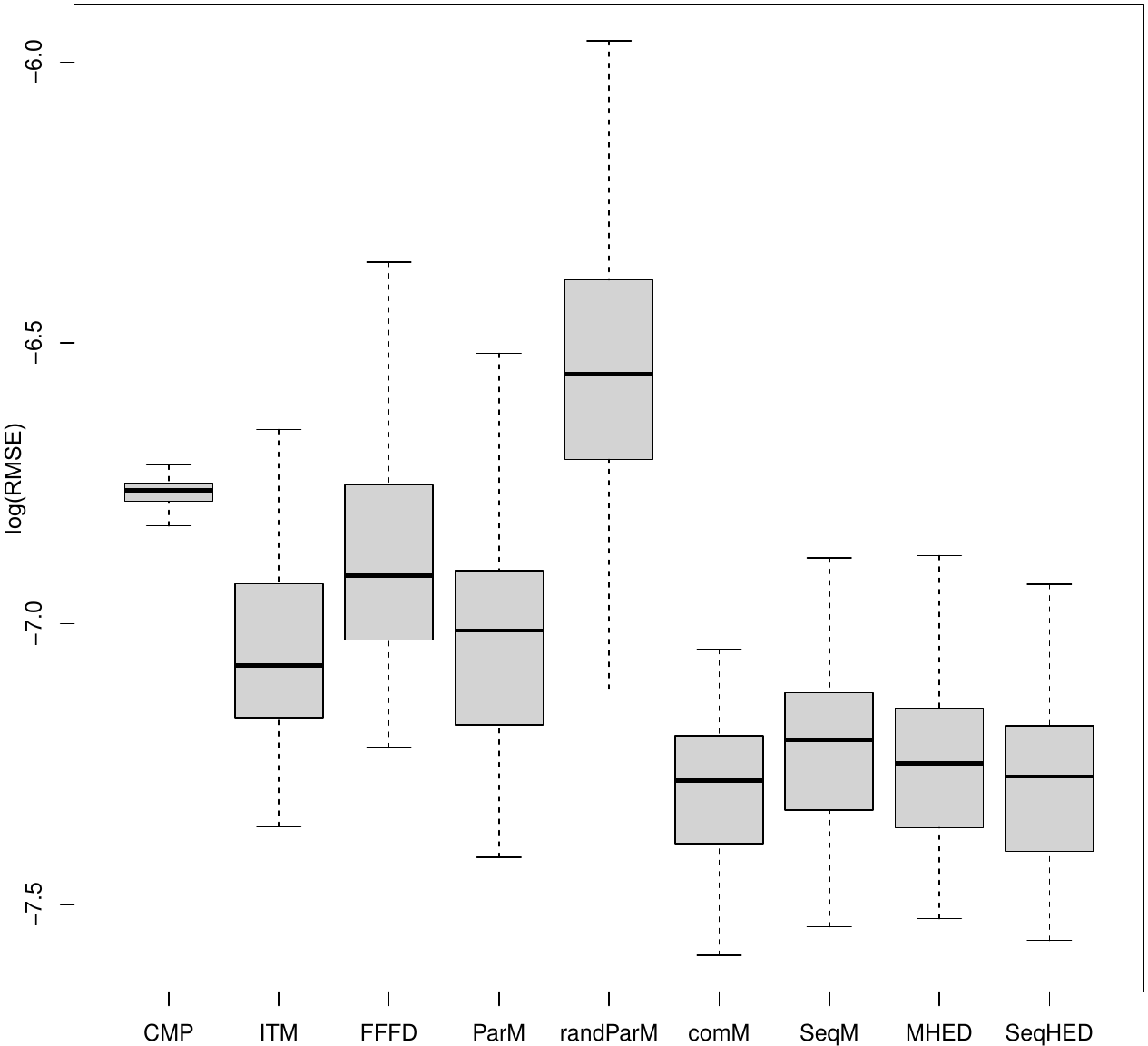}
        \caption{EzGP}
  \end{subfigure}
  \begin{subfigure}[b]{0.32\textwidth}
        \includegraphics[width=\textwidth]{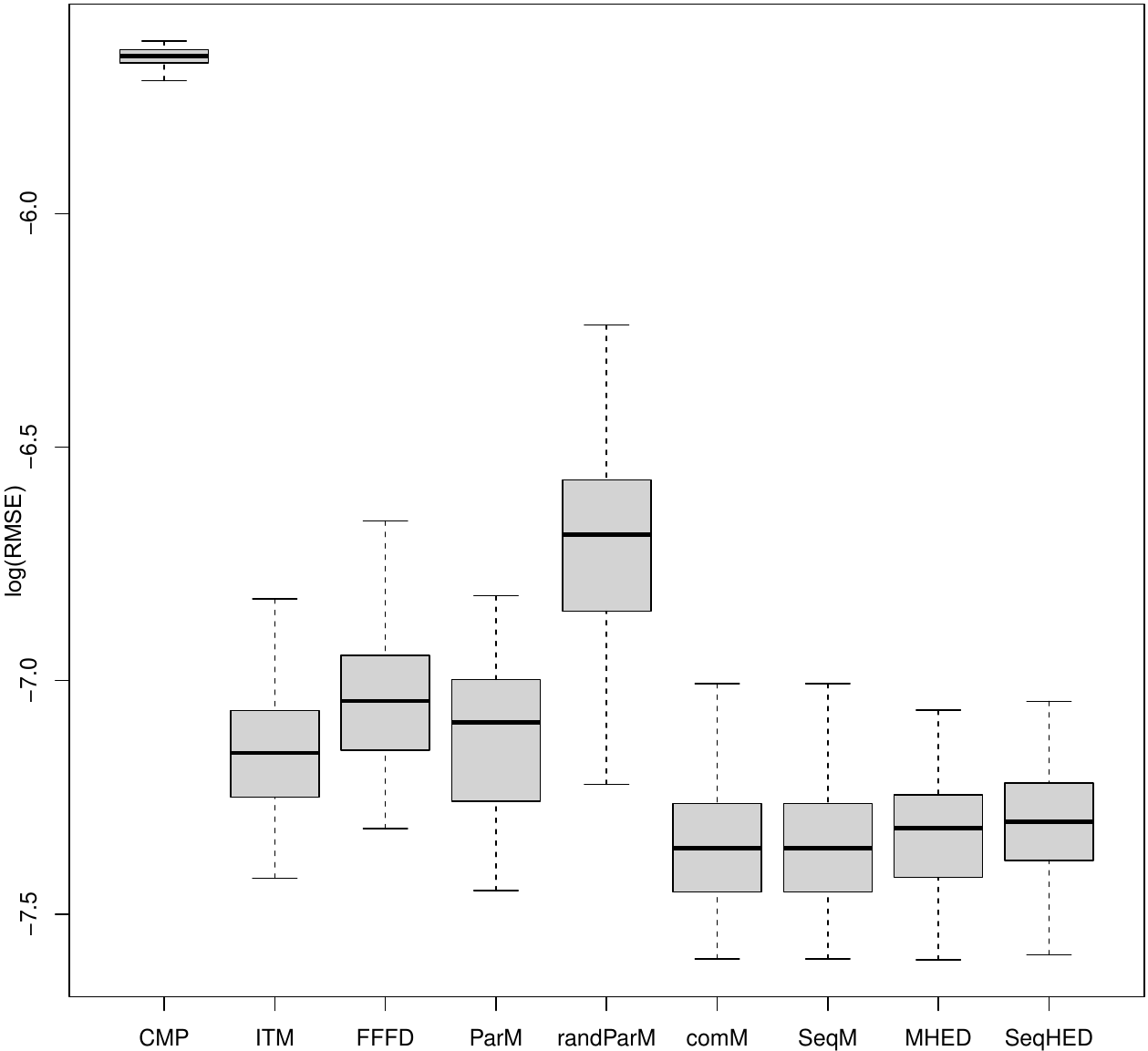}
        \caption{EEzGP}
  \end{subfigure}
 \caption{\ Comparison of cell survival rate prediction errors based on different mixture design methods and modeling methods.}
\label{fig5}
\end{figure}

We generate 200 uniform test samples across different levels of the total dose of the drugs to evaluate the performance of each design scheme, so the total test sample size is $M=1400$. If the total dose $ C $ is treated as a continuous variable, we build an ordinary Gaussian process (OGP) model with constant mean using $ (\theta_1, \theta_2, C) $ as inputs, and compute the total RMSE defined in \eqref{RMSEcri}. Each scenario is repeated 50 times, and the results are shown in Figure\ref{fig5}(a). From Figure\ref{fig5}(a), we observe that the sliced space-filling design with mixtures proposed in this paper achieve smaller errors than other methods, with MHED performing the best. Although randParM and ParM exhibit similar uniformity in the full design, ParM outperforms randParM significantly, especially when there are large differences in response surfaces.  

If $ C $ is treated as a qualitative variable, it cannot be directly combined with mixture components for analysis using a ordinary Gaussian process model. Unlike Example \ref{ex2}, when response surfaces of mixtures differ, we cannot merge mixture component data for modeling. \cite{EzGPpaper} proposed the EzGP and EEzGP models for modeling and predicting experimental data, which are suitable for problems involving both quantitative and qualitative factors. The EzGP model assigns different kernel parameters to different levels of qualitative variables to distinguish response surface differences, while the EEzGP model uses the same parameters to simplify the model and assumes similarity across levels of the same qualitative factor (for details, refer to \citealp{EzGPpaper} and the R package \texttt{EzGP} \citealp{EzGPr}).  

We employ EzGP and EEzGP to model all experimental data, computing prediction errors on the test set (excluding errors from subdesign-based modeling). Results in Figure~\ref{fig5}(b-c) show that the four energy distance-based criteria outperform other methods while exhibiting comparable performance. For EzGP, the comM design performs best, as it accurately estimates normal cell mortality across dose levels; the other three methods also perform well due to good subdesign uniformity. For EEzGP, designs with good full uniformity show improved prediction accuracy. In contrast, the CMP method exhibits poor full uniformity, leading to reduced prediction accuracy. The CDM method is omitted from figures owing to its large variance and inferior performance.

\section{Summary and future research directions}

Leveraging energy distance decomposition theory, this paper develops three methods for constructing representative point sets with partitioned structures. We establish that under specific conditions, both the complete point sets and their partitioned subsets generated by these methods converge provably to the target distribution. Furthermore, within the hybrid energy distance criterion framework and Majorization-Minimization optimization, we present multiple construction algorithms for mixture partition designs with rigorous convergence guarantees.

Numerical experiments demonstrate that the proposed sliced space-filling design with mixtures outperform existing methods across multiple uniformity criteria. Under diverse realistic model Lemmas, the two hybrid energy distance-based methods achieve the smallest mean squared prediction error. Crucially, these methods offer two significant advantages:
(i) applicability to experimental regions defined by convex sets;
(ii) flexible specification of run sizes and target distributions for sub-designs, greatly enhancing practical utility.

As the number of process variables and their levels increases, the required experimental points for mixture partition designs grow exponentially with the number of level combinations. To mitigate this scaling issue, we may develop novel design approaches by integrating the proposed energy distance decomposition theory with Marginal Coupled Design (MCD) \citep{MCD} and Double Coupled Design (DCD) \citep{DCD} methodologies. These integrated methods preserve high uniformity in mixture designs across:
(i) different levels of individual process variables, and
(ii) level combinations of any two process variables.
In addition, the isotropic nature of Euclidean distance may lead to suboptimal distributions of representative points in low-dimensional projection spaces when used in energy distance criteria. Consequently, developing appropriate metrics to enhance projection-space uniformity remains an open research problem. Investigating these theoretical aspects and designing corresponding optimization algorithms present promising directions for future work.

\section*{Disclosure statement}
No potential conflict of interest was reported by the authors.

\section*{Data availability statement}
The authors confirm that the data supporting the findings of this study are available within the article and source codes are available on GitHub (\url{https://github.com/doenew/SlicedSFDwithmixtures}).

\phantomsection\label{supplementary-material}
\bigskip

\begin{center}

{\large\bf SUPPLEMENTARY MATERIAL}

\end{center}

\begin{description}
\item[Proofs of technical results:] The proofs of Theorems 1-8 and Propositions 2.
\item[Visualizations of sliced space-filling design with mixtures:]
The ternary scatter plots of designs generated by seven methods in example \ref{ex1}.
\end{description}

 \bibliography{bibliography.bib}

\end{document}


\maketitle
\thispagestyle{empty} 


\section*{Additional proofs of technical results}
\label{appA}

\noindent\textbf{Proof of Theorem 1 }
Let $ \mathcal{I}_{n_k} $ denote the index set corresponding to the sample points in the full point set $ \mathcal{P}_n $ that are contained in the subset $ \mathcal{P}_{n_k} $. 
According to Definition 2 and (3),
\begin{equation*}
\begin{aligned}
E(F,F_{\mathcal{P}_n})=&\frac{2}{n}\sum_{i=1}^n\mathbb{E}[\|\boldsymbol{x}_i-Y\|_2]-\frac{1}{n^2}
  \sum_{i=1}^n\sum_{j=1}^n\|\boldsymbol{x}_i-\boldsymbol{x}_j\|_2-\mathbb{E}\|Y-Y'\|_2\\
  =&\frac{2}{n}\sum_{k=1}^K\sum_{i\in\mathcal{I}_{n_k}}\mathbb{E}[\|\boldsymbol{x}_i-Y\|_2]-\frac{1}{n^2}\sum_{k=1}^K\sum_{i,j\in\mathcal{I}_{n_k}}\|\boldsymbol{x}_i-\boldsymbol{x}_j\|_2\\
& -\frac{2}{n^2}\sum_{1\leq {k_1}<{k_2}\leq K}\sum_{i\in\mathcal{I}_{n_{k_1}}}\sum_{j\in\mathcal{I}_{n_{k_2}}}
\|\boldsymbol{x}_i-\boldsymbol{x}_j\|_2-\frac{1}{n}\sum_{k=1}^Kn_k\mathbb{E}\|Y-Y'\|_2\\
=&\sum_{k=1}^K\frac{n_k}{n}\left\{\frac{2}{n_k}\sum_{i\in\mathcal{I}_{n_k}}\mathbb{E}[\|\boldsymbol{x}_i-Y\|_2]-\frac{1}{n_k^2}\sum_{i,j\in\mathcal{I}_{n_k}}
\|\boldsymbol{x}_i-\boldsymbol{x}_j\|_2-\mathbb{E}\|Y-Y'\|_2\right\}\\
&-\sum_{1\leq {k_1}<{k_2}\leq K}\frac{n_{k_1}n_{k_2}}{n^2}\left\{\frac{2}{n_{k_1}n_{k_2}}\sum_{i\in\mathcal{I}_{n_{k_1}}}\sum_{j\in\mathcal{I}_{n_{k_2}}}
\|\boldsymbol{x}_i-\boldsymbol{x}_j\|_2\right.\\
&\left.-\frac{1}{n_{k_1}^2}\sum_{i,j\in\mathcal{I}_{n_{k_1}}}
\|\boldsymbol{x}_i-\boldsymbol{x}_j\|_2-\frac{1}{n_{k_2}^2}\sum_{i,j\in\mathcal{I}_{n_{k_2}}}
\|\boldsymbol{x}_i-\boldsymbol{x}_j\|_2\right\}\\
=&\sum_{k=1}^K{\frac{n_k}{n}}E\left( F,F_{\mathcal{P}_{n_k}} \right) -\sum_{1\leq {k_1}<{k_2}\leq K}{\frac{n_{k_1}n_{k_2}}{n^2}}E\left( F_{\mathcal{P}_{n_{k_1}}},F_{\mathcal{P}_{n_{k_2}}} \right)\nonumber \\
=&\sum_{k=1}^K{\frac{n_k}{n}}E\left( F,F_{\mathcal{P}_{n_k}} \right) -\sum_{{k_1}=1}^K{\sum_{{k_2}=1}^K{\frac{n_{k_1}n_{k_2}}{2n^2}}}E\left( F_{\mathcal{P}_{n_{k_1}}},F_{\mathcal{P}_{n_{k_2}}} \right).
\end{aligned}
\end{equation*}
To establish the third equality above, it suffices to verify that the coefficients preceding $ \sum_{i,j \in \mathcal{I}_{n_k}} \delta(\boldsymbol{x}_i, \boldsymbol{x}_j) $ are equal. The coefficient on the left-hand side of the equation is $ -\frac{1}{n^2} $, while the coefficient on the right-hand side is:
\begin{equation*}
    -\frac{1}{n n_k} + \sum_{1 \leq k_1 < k} \frac{1}{n^2} \frac{n_{k_1}}{n_k} + \sum_{k < k_2 \leq K} \frac{1}{n^2} \frac{n_{k_2}}{n_k} = -\frac{1}{n n_k} + \frac{1}{n^2 n_k}(n - n_k) = -\frac{1}{n^2}.
\end{equation*}

Let $ F = F_{\mathcal{P}_n} $ in above equality. By Lemma 1, $ E(F_{\mathcal{P}_n}, F_{\mathcal{P}_n})$ $= 0 $, then
\begin{equation*}
    \sum_{{k_1}=1}^K{\sum_{{k_2}=1}^K{\frac{n_{k_1}n_{k_2}}{2n^2}}}E\left( F_{\mathcal{P}_{n_{k_1}}},F_{\mathcal{P}_{n_{k_2}}} \right)=\sum_{k=1}^K{\frac{n_k}{n}}E\left( F_{\mathcal{P}_n},F_{\mathcal{P}_{n_k}} \right), 
\end{equation*}
Therefore, the last equality in (6) holds.
\qed

\vspace{0.5cm}

\noindent\textbf{Proof of Theorem 2 }
By definition, $ c_k \geq 0 $ and $ \sum_{k=1}^K c_k = 1 $. According to Proposition 1 and the definition of the set $ S $:
\begin{equation*}
    \begin{aligned}
      0 \leq E\left(F, F_{\mathcal{P}_n^*}\right) &\leq \sum_{k=1}^K \frac{n_k(n)}{n} E\left(F, F_{\mathcal{P}_{n_k(n)}^*}\right) \\ &= \sum_{k \in S} \frac{n_k(n)}{n} E\left(F, F_{\mathcal{P}_{n_k(n)}^*}\right) + \sum_{k \notin S} \frac{n_k(n)}{n} E\left(F, F_{\mathcal{P}_{n_k(n)}^*}\right).  
    \end{aligned}
\end{equation*}

According to Lemma 1 and Lemma 2 $\forall\ k \in S $, we have $ F_{\mathcal{P}_{n_k}^*} \stackrel{d}{\longrightarrow} F $ and $ \lim_{n \to \infty} E\left(F, F_{\mathcal{P}_{n_k(n)}^*}\right)$ $ = 0 $ for all $ k \in S $. For every $ k \notin S $, $ \lim_{n \to \infty} \frac{n_k(n)}{n} = 0 $ and $ E\left(F, F_{\mathcal{P}_{n_k(n)}^*}\right) < \infty $. Thus:
\[
\lim_{n \to \infty} \sum_{k \in S} \frac{n_k(n)}{n} E\left(F, F_{\mathcal{P}_{n_k(n)}^*}\right) = 0.
\]
Therefore, $ \lim_{n \to \infty} E\left(F, F_{\mathcal{P}_n^*}\right) = 0 $, and $ \lim_{n\rightarrow\infty}F_{\mathcal{P}_n^*} = F $ by Lemma 2.  \qed

\vspace{0.5cm}

\noindent\textbf{Proof of Proposition 2 }
When $ \mathcal{P}_{n_1} = \cdots = \mathcal{P}_{n_K} $, the sample sizes and empirical distribution functions of all subsets are identical. Substituting into (6) yields the conclusion stated in the proposition.  
\qed

\vspace{0.5cm}

\begin{slemma}\label{sLemma1}
Let $F_{\mathcal{P}_n}$ be the empirical distribution function of sample set ${\mathcal{P}_n}$ and $\lim_{n\rightarrow\infty}F_{\mathcal{P}_n}=F$. Suppose $ n = \sum_{k=1}^K n_k(n) $ and $ \lim_{n \to \infty} \frac{n_k(n)}{n} = c_k > 0 $ for all $ k = 1, \ldots, K $. Randomly partition $\mathcal{P}_n$ into $K$ subsets $\mathcal{P}_n = \bigcup_{k = 1}^K \tilde{\mathcal{P}}_{n_k(n)}$, then $\lim_{n\rightarrow\infty}F_{\tilde{\mathcal{P}}_{n_k(n)}}=F$, for all $ k = 1, \ldots, K $.
\end{slemma}

\noindent\textbf{Proof of Lemma \ref{sLemma1} }
Let $ \{X_1, \dots, X_n\} $ denote the sample points in $ \mathcal{P}_n $. Define assignment indicators $ Z_{i,k} = \mathbf{1}_{\{X_i \in \tilde{\mathcal{P}}_{n_k(n)}\}} $ for $ i = 1,\dots,n $ and $ k = 1,\dots,K $. The empirical distribution of subset $ \tilde{\mathcal{P}}_{n_k(n)} $ is:
\begin{equation*}
  F_{\tilde{\mathcal{P}}_{n_k(n)}}(\boldsymbol{x}) = \frac{1}{n_k(n)} \sum_{i=1}^n Z_{i,k} \mathbf{1}_{\{X_i \leq \boldsymbol{x}\}}.  
\end{equation*}

Conditional on $ \mathcal{P}_n $:
\begin{equation*}
\begin{split}
       \mathbb{E}\left[ F_{\tilde{\mathcal{P}}_{n_k(n)}}(\boldsymbol{x}) \mid \mathcal{P}_n \right] &= \frac{1}{n_k(n)} \sum_{i=1}^n \mathbb{E}[Z_{i,k} \mid \mathcal{P}_n] \mathbf{1}_{\{X_i \leq \boldsymbol{x}\}} \\
       &= \frac{1}{n_k(n)} \sum_{i=1}^n \frac{n_k(n)}{n} \mathbf{1}_{\{X_i \leq \boldsymbol{x}\}} = F_{\mathcal{P}_n}(\boldsymbol{x}),  
\end{split}
\end{equation*}
since $ \mathbb{E}[Z_{i,k} \mid \mathcal{P}_n] = \frac{n_k(n)}{n} $ by symmetry of random partitioning.
\begin{equation*}
        \mathrm{Var}\left( F_{\tilde{\mathcal{P}}_{n_k(n)}}(\boldsymbol{x}) \mid \mathcal{P}_n \right) = \frac{1}{n_k(n)^2} \mathrm{Var}\left( \sum_{i=1}^n Z_{i,k} \mathbf{1}_{\{X_i \leq \boldsymbol{x}\}} \mid \mathcal{P}_n \right).
\end{equation*}
Using hypergeometric variance properties:
 \begin{equation*}
   \mathrm{Var}\left( \sum_{i=1}^n Z_{i,k} Y_i \mid \mathcal{P}_n \right) = \frac{n_k(n)(n - n_k(n))}{n(n-1)} M(1 - M/n),     
   \end{equation*}
where $ M = n F_{\mathcal{P}_n}(\boldsymbol{x}) $ and $ Y_i = \mathbf{1}_{\{X_i \leq \boldsymbol{x}\}} $. Thus:
\begin{equation*}
        \mathrm{Var}\left( F_{\tilde{\mathcal{P}}_{n_k(n)}}(\boldsymbol{x}) \mid \mathcal{P}_n \right) \leq \frac{1}{4n_k(n)} \left(1 - \frac{n_k(n)}{n}\right) \frac{n}{n-1}.
\end{equation*}
As $ n \to \infty $:
\begin{enumerate}
    \item[(i)] $ \frac{n_k(n)}{n} \to c_k > 0 $ implies $ n_k(n) \to \infty $;
    \item[(ii)] $ F_{\mathcal{P}_n}(\boldsymbol{x}) \to F(\boldsymbol{x}) $ at continuity points $\boldsymbol{x} $ of $ F $;
    \item[(iii)] Conditional variance $\mathrm{Var}\left( F_{\tilde{\mathcal{P}}_{n_k(n)}}(\boldsymbol{x}) \mid \mathcal{P}_n \right) \to 0 $.
\end{enumerate}
By Chebyshev's inequality, for any $ \epsilon > 0 $:
\[
\mathbb{\mathcal{P}}\left( \left| F_{\tilde{\mathcal{P}}_{n_k(n)}}(\boldsymbol{x}) - F_{\mathcal{P}_n}(\boldsymbol{x}) \right| \geq \epsilon \mid \mathcal{P}_n \right) \leq \frac{\text{Var}\left( F_{\tilde{\mathcal{P}}_{n_k(n)}}(\boldsymbol{x}) \mid \mathcal{P}_n \right)}{\epsilon^2} \to 0.
\]
Hence $ F_{\tilde{\mathcal{P}}_{n_k(n)}}(\boldsymbol{x}) \xrightarrow{{P}} F_{\mathcal{P}_n}(\boldsymbol{x}) $, and since $ F_{\mathcal{P}_n}(\boldsymbol{x}) \to F(\boldsymbol{x}) $, we have $ F_{\tilde{\mathcal{P}}_{n_k(n)}}(\boldsymbol{x}) \xrightarrow{{P}} F(\boldsymbol{x}) $ at continuity points of $ F $ by the Slutsky's theorem. This implies $ F_{\tilde{\mathcal{P}}_{n_k(n)}} \xrightarrow{d} F $ by the Portmanteau theorem \citep{Probability}. \qed

\vspace{0.5cm}

\noindent\textbf{Proof of Theorem 3 }
By Lemma 1 and Lemma 2, we have $ \lim_{n\rightarrow\infty}F_{\mathcal{P}_n^*} = F $ and $ \lim_{n \to \infty} E\left(F, F_{\mathcal{P}_n^*}\right) = 0 $. Let $\tilde{\mathcal{P}}_{n_1(n)},\dots,\tilde{\mathcal{P}}_{n_K(n)}$ be a random partition of $\mathcal{P}_n^*$, then the empirical distribution function $F_{\tilde{\mathcal{P}}_{n_K(n)}}\stackrel{d}{\longrightarrow} F$ of $\tilde{\mathcal{P}}_{n_K(n)}$ by Lemma \ref{sLemma1}, $k=1,\dots,K$. Let $\mathcal{P}^*_{n_1(n)},\dots,\mathcal{P}^*_{n_K(n)}$ be a solution of (9), we have

Thus, the objective function in (9) satisfies:  
\begin{equation*}
  0\leq\lim_{n \to \infty} \sum_{k=1}^K \frac{n_k(n)}{n} E\left(F, F_{\mathcal{P}_{n_k(n)}^*}\right) \leq 
  \lim_{n \to \infty} \sum_{k=1}^K \frac{n_k(n)}{n} E\left(F, F_{\tilde{\mathcal{P}}_{n_k(n)}}\right)=0.  
\end{equation*}
Therefore $ \lim_{n_k(n) \to \infty} E\left(F, F_{\mathcal{P}_{n_k(n)}^*}\right) = 0 $, $ \forall\ k = 1, \ldots, K $.   \qed

\vspace{0.5cm}

\noindent\textbf{Proof of Theorem 4 }
By Proposition 1, we have:  
\begin{equation*}
 E^\lambda\left(F, F_{\mathcal{P}_n}\right) \leq \sum_{k=1}^K \frac{n_k(n)}{n} E\left(F, F_{\mathcal{P}_{n_k(n)}}\right).   
\end{equation*}
Define:  
\begin{equation*}
 \tilde{\mathcal{P}}_n=\underset{\mathcal{P}_{n_1(n)},\dots,\mathcal{P}_{n_K(n)} \subseteq \mathcal{X}}{\textup{argmin}} \sum_{k=1}^K{\frac{n_k(n)}{n}}E\left( F,F_{\mathcal{P}_{n_k(n)}} \right),
\end{equation*}
where $ \widetilde{\mathcal{P}}_n = \bigcup_{k=1}^K \widetilde{\mathcal{P}}_{n_k(n)} $.  
Based on the separability of $ \sum_{k=1}^K \frac{n_k}{n} E \left( F, F_{\mathcal{P}_{n_k}} \right) $ and Lemma 2, we derive:  
\begin{equation*}
  0\leq \lim_{n\rightarrow \infty} E^{\lambda}\left(F,F_{\mathcal{P}^*_n}\right)\leq \lim_{n\rightarrow \infty}\sum_{k=1}^K{\frac{n_k}{n}}E\left( F,F_{\tilde{\mathcal{P}}_{n_k(n)}}\right)=\sum_{k=1}^Kc_k\lim_{n\rightarrow \infty}E\left( F,F_{\tilde{\mathcal{P}}_{n_k(n)}}\right)=0.
\end{equation*}
Thus, 
\begin{equation*}
    \lim_{n \to \infty} E^\lambda \left( F, F_{\mathcal{P}_n^*} \right) = \lim_{n \to \infty} \lambda E(F, F_{\mathcal{P}_n^*}) + (1-\lambda) \sum_{k=1}^K \frac{n_k}{n} E \left( F, F_{\mathcal{P}_{n_k}^*} \right) = 0.
\end{equation*}
From the non-negativity of the energy distance, it follows that:
\begin{equation*}
  \lim_{n \to \infty} E(F, F_{\mathcal{P}_n^*}) = 0, \quad \lim_{n \to \infty} E \left( F, F_{\mathcal{P}_{n_k}^*} \right) = 0 \quad (k = 1, \ldots, K).  
\end{equation*}
By Lemma 2, we conclude $ \lim_{n\rightarrow\infty}F_{\mathcal{P}_n^*} = F $ and $\lim_{n\rightarrow\infty} F_{\mathcal{P}_{n_k(n)}^*}  = F $ for all $ k = 1, \ldots, K $.  \qed

\vspace{0.5cm}

\noindent\textbf{Proof of Theorem 5 }
By the arithmetic-geometric mean inequality, 
\begin{equation*}
  \frac{\|\boldsymbol{x}_i-\boldsymbol{y}_m
\|_2^2}{\|\boldsymbol{x}^{(t)}_i-\boldsymbol{y}_m\|_2}+\|\boldsymbol{x}^{(t)}_i-\boldsymbol{y}_m\|_2\geq 2\sqrt{\frac{\|\boldsymbol{x}_i-\boldsymbol{y}_m
\|_2^2}{\|\boldsymbol{x}^{(t)}_i-\boldsymbol{y}_m\|_2}\|\boldsymbol{x}^{(t)}_i-\boldsymbol{y}_m\|_2}=2\|\boldsymbol{x}_i-\boldsymbol{y}_m
\|_2.
\end{equation*}
Let $f(\boldsymbol{a},\boldsymbol{b})=\|\boldsymbol{a}-\boldsymbol{b}\|_2$, $\forall \theta\in[0,1]$,
\begin{equation*}
\begin{aligned}
  f(\theta \boldsymbol{a}_1+(1-\theta)\boldsymbol{a}_2,\theta \boldsymbol{b}_1+(1-\theta)\boldsymbol{b}_2)=&\|\theta \boldsymbol{a}_1+(1-\theta)\boldsymbol{a}_2-\theta \boldsymbol{b}_1+(1-\theta)\boldsymbol{b}_2\|_2\\
  =& \|\theta (\boldsymbol{a}_1-\boldsymbol{b}_1)+(1-\theta)(\boldsymbol{a}_2-\boldsymbol{b}_2)\|_2\\
  \leq & \theta \|\boldsymbol{a}_1-\boldsymbol{b}_1\|_2+(1-\theta)\|\boldsymbol{a}_2-\boldsymbol{b}_2\|_2\\
  =&\theta f(\boldsymbol{a}_1,\boldsymbol{b}_1)+(1-\theta)f(\boldsymbol{a}_2,\boldsymbol{b}_2).
\end{aligned}
\end{equation*} 
Thus, $ f(\boldsymbol{a}, \boldsymbol{b}) = \|\boldsymbol{a} - \boldsymbol{b}\|_2 $ is a convex function. Performing a first-order Taylor expansion of $ \|\boldsymbol{x}_i - \boldsymbol{x}_j\|_2 $ around $ (\boldsymbol{x}_i^{(t)}, \boldsymbol{x}_j^{(t)}) $, we have:  
\begin{equation*}
  \|\boldsymbol{x}_i - \boldsymbol{x}_j\|_2 \geq \|\boldsymbol{x}_i^{(t)} - \boldsymbol{x}_j^{(t)}\|_2 + \frac{(\boldsymbol{x}_i^{(t)} - \boldsymbol{x}_j^{(t)})^T (\boldsymbol{x}_i - \boldsymbol{x}_i^{(t)})}{\|\boldsymbol{x}_i^{(t)} - \boldsymbol{x}_j^{(t)}\|_2} + \frac{(\boldsymbol{x}_j^{(t)} - \boldsymbol{x}_i^{(t)})^T (\boldsymbol{x}_j - \boldsymbol{x}_j^{(t)})}{\|\boldsymbol{x}_i^{(t)} - \boldsymbol{x}_j^{(t)}\|_2}.  
\end{equation*}
Applying the above inequality to $ g(\mathcal{P}_n \mid \mathcal{P}_n^{(t)}) $ and using the symmetry of the index sets in the summations, we derive $ g(\mathcal{P}_n \mid \mathcal{P}_n^{(t)}) \geq h(\mathcal{P}_n) $ for all $ \mathcal{P}_n \subseteq \mathcal{X} $. Setting $ \mathcal{P}_n = \mathcal{P}_n^{(t)} $ in $ g(\mathcal{P}_n \mid \mathcal{P}_n^{(t)}) $, it is straightforward to verify $ g(\mathcal{P}_n^{(t)} \mid \mathcal{P}_n^{(t)}) = h(\mathcal{P}_n^{(t)}) $. Thus, $ g(\mathcal{P}_n \mid \mathcal{P}_n^{(t)}) $ is a surrogate function of $ h(\mathcal{P}_n) $ at the current iteration $ \mathcal{P}_n^{(t)} $.  

Since $ g(\mathcal{P}_n \mid \mathcal{P}_n^{(t)}) $ is a quadratic function in all optimization variables $ \boldsymbol{x}_{\mathcal{P}_{n_k}, i} $ for $ i = 1, \cdots, n_k $ and $ k = 1, \ldots, K $, we take the partial derivative with respect to each $ \boldsymbol{x}_{\mathcal{P}_{n_k}, i} $ and set it to zero:  
\begin{equation*}
   \begin{split}
\frac{\partial g(\mathcal{P}_n \mid \mathcal{P}_n^{(t)})}{\partial \boldsymbol{x}_{\mathcal{P}_{n_k}, i}} &= \frac{1}{nN} \sum_{m=1}^N \frac{2(\boldsymbol{x}_{\mathcal{P}_{n_k}, i} - \boldsymbol{y}_m)}{\|\boldsymbol{x}_{\mathcal{P}_{n_k}, i}^{(t)} - \boldsymbol{y}_m\|_2} - \left( \frac{1-\lambda}{nn_k} + \frac{\lambda}{n^2} \right) \sum_{\substack{j \in \mathcal{I}_{n_k} \\ j \neq i}} \frac{2(\boldsymbol{x}_{\mathcal{P}_{n_k}, i}^{(t)} - \boldsymbol{x}_j^{(t)})}{\|\boldsymbol{x}_{\mathcal{P}_{n_k}, i}^{(t)} - \boldsymbol{x}_j^{(t)}\|_2} \\
&\quad - \frac{\lambda}{n^2} \sum_{\substack{k_2=1 \\ k_2 \neq k}}^K \sum_{j \in \mathcal{I}_{n_{k_2}}} \frac{2(\boldsymbol{x}_{\mathcal{P}_{n_k}, i}^{(t)} - \boldsymbol{x}_j^{(t)})}{\|\boldsymbol{x}_{\mathcal{P}_{n_k}, i}^{(t)} - \boldsymbol{x}_j^{(t)}\|_2} = 0.
\end{split} 
\end{equation*}
Simplifying the above equation yields the result in (17).   \qed

\vspace{0.5cm}

\noindent\textbf{Proof of Theorem 6}
By the convergence conditions of the MM algorithm \citep{Mak2018,MMcon} and the continuity of the surrogate function $ g(\mathcal{P}_n \mid \mathcal{P}_n^{(t)}) $, to prove the conclusion in Theorem 6, it suffices to show that the directional derivatives of $ g(\mathcal{P}_n \mid \mathcal{P}_n^{(t)}) $ and $ h(\mathcal{P}_n) $ coincide in any feasible direction at $ \mathcal{P}_n^{(t)} $. We verify this as follows.  
For any feasible direction set:  
\begin{equation*}
  \left\{ \boldsymbol{d} \mid \boldsymbol{d} \subseteq \mathbb{R}^{np} \text{ and } \boldsymbol{x}_{\mathcal{P}_{n_k}, i} + \boldsymbol{d}_{\mathcal{P}_{n_k}, i} \in \mathcal{X} \text{ for all } i = 1, \ldots, n_k \text{ and } k = 1, \ldots, K \right\},  
\end{equation*}
we compute the directional derivative:  
\begin{equation*}
\begin{split}
    &\left. g'\left( \mathcal{P}_n \mid \mathcal{P}_n^{(t)}; \boldsymbol{d} \right) \right|_{\mathcal{P}_n = \mathcal{P}_n^{(t)}} \\
&= \lim_{\lambda \to 0^+} \frac{g\left( \{\boldsymbol{x}_{\mathcal{P}_{n_k}, i} + \lambda \boldsymbol{d}_{\mathcal{P}_{n_k}, i}\}_{i,k} \mid \mathcal{P}_n^{(t)} \right) - g\left( \{\boldsymbol{x}_{\mathcal{P}_{n_k}, i}\}_{i,k} \mid \mathcal{P}_n^{(t)} \right)}{\lambda} \\
&= \frac{2}{nN} \sum_{m=1}^N \sum_{k=1}^K \sum_{i \in \mathcal{I}_{n_k}} \left( \frac{(\boldsymbol{x}_i^{(t)} - \boldsymbol{y}_m)^T \boldsymbol{d}_i}{\|\boldsymbol{x}_i^{(t)} - \boldsymbol{y}_m\|_2} \right) \\
&\quad - \sum_{k=1}^K \left( \frac{1-\lambda}{nn_k} + \frac{\lambda}{n^2} \right) \sum_{\substack{i,j \in \mathcal{I}_{n_k} \\ i \neq j}} \left[ \frac{2(\boldsymbol{x}_i^{(t)} - \boldsymbol{x}_j^{(t)})^T \boldsymbol{d}_i}{\|\boldsymbol{x}_i^{(t)} - \boldsymbol{x}_j^{(t)}\|_2} \right] \\
&\quad - \frac{\lambda}{n^2} \sum_{k_1=1}^K \sum_{\substack{k_2=1 \\ k_2 \neq k_1}}^K \sum_{i \in \mathcal{I}_{n_{k_1}}} \sum_{j \in \mathcal{I}_{n_{k_2}}} \left[ \frac{2(\boldsymbol{x}_i^{(t)} - \boldsymbol{x}_j^{(t)})^T  \boldsymbol{d}_i}{\|\boldsymbol{x}_i^{(t)} - \boldsymbol{x}_j^{(t)}\|_2} \right] \\
&= \left. h'\left( \mathcal{P}_n; \boldsymbol{d} \right) \right|_{\mathcal{P}_n = \mathcal{P}_n^{(t)}}.
\end{split}
\end{equation*}

Therefore, the limit point set of $ \{\mathcal{P}_n^{(t)}\}_{t=0}^\infty $ converges to a stationary solution of (15).   \qed

\vspace{0.5cm}

\noindent\textbf{Proof of Theorem 7 and Theorem 8}
The proofs of these two theorems are respectively similar to those of Theorem 5 and Theorem 6, and are omitted here.

\section*{Additional visualizations of sliced space-filling design with mixtures} \label{appB}

The scatter plots of sliced space-filling design with mixtures in $T^3$ generated by seven methods are shown in Fig.\ref{figITM} to Fig.\ref{figSeqHED}. Different shapes or colors correspond to different subdesigns.

\setcounter{figure}{0}

\begin{figure}[htbp]
    \centering
    \begin{subfigure}[b]{0.49\textwidth}
        \includegraphics[width=\textwidth]{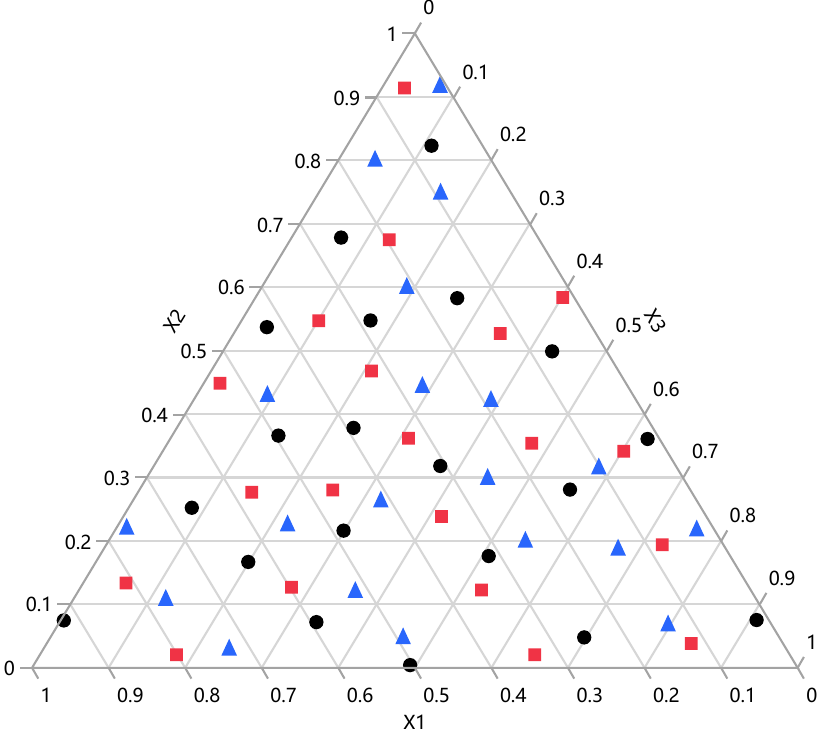}
        \caption{full design}
     \end{subfigure}
    \begin{subfigure}[b]{0.49\textwidth}
        \includegraphics[width=\textwidth]{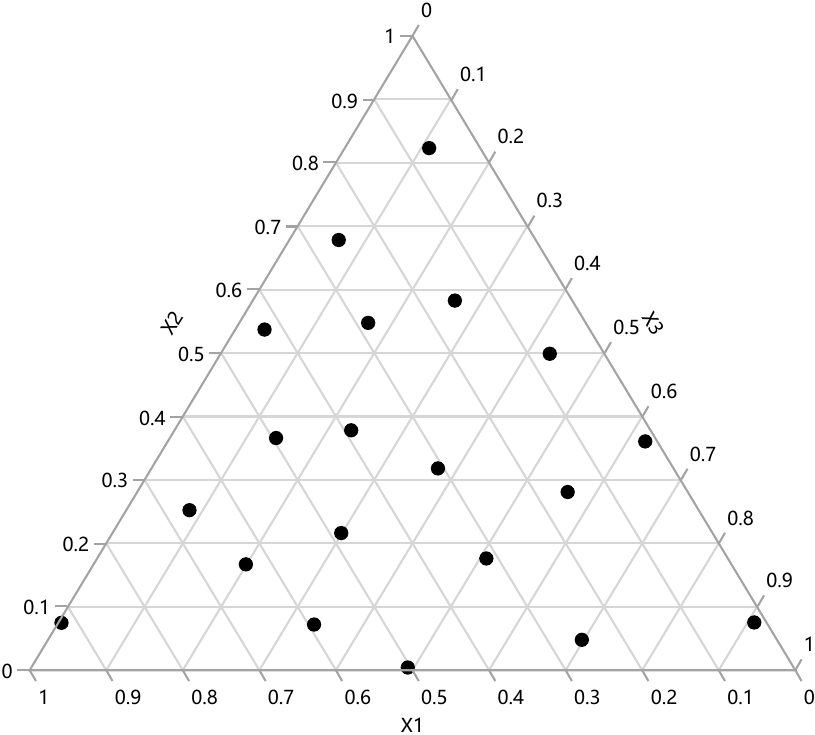}
        \caption{ slice 1}
  \end{subfigure}
  \begin{subfigure}[b]{0.49\textwidth}
        \includegraphics[width=\textwidth]{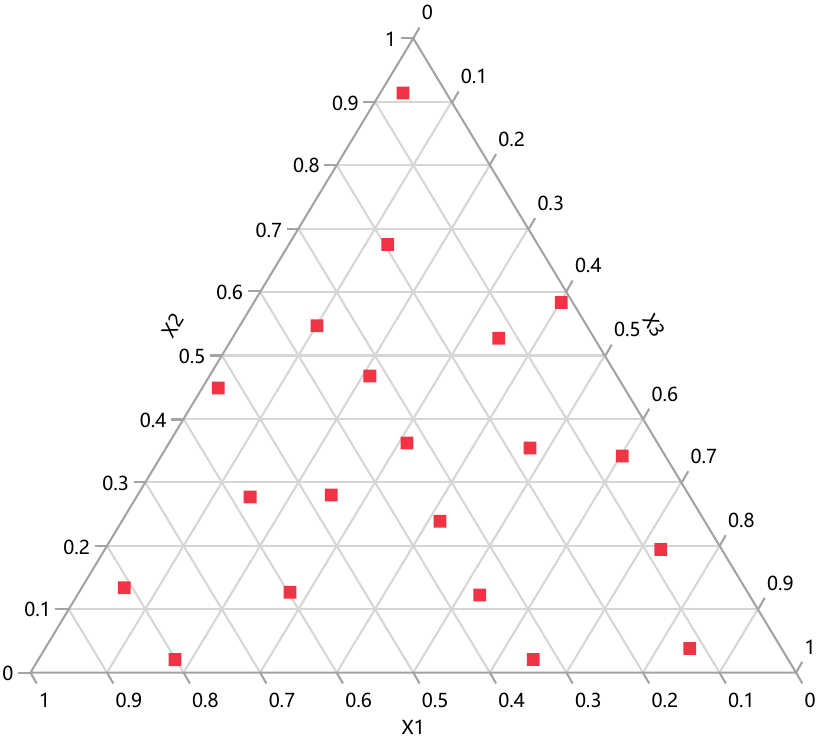}
        \caption{slice 2}
  \end{subfigure}
    \begin{subfigure}[b]{0.49\textwidth}
        \includegraphics[width=\textwidth]{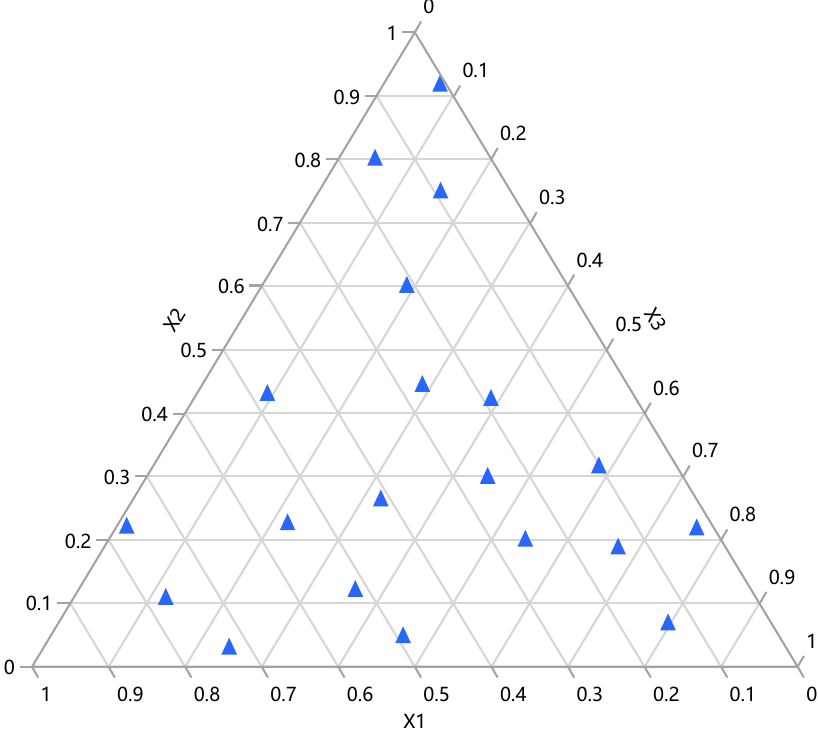}
        \caption{slice 3}
  \end{subfigure}
 \caption{\ Scatter plot of sliced space-filling design with mixtures generated by ITM}
\label{figITM}
\end{figure}
\begin{figure}[htbp]
    \centering
    \begin{subfigure}[b]{0.49\textwidth}
        \includegraphics[width=\textwidth]{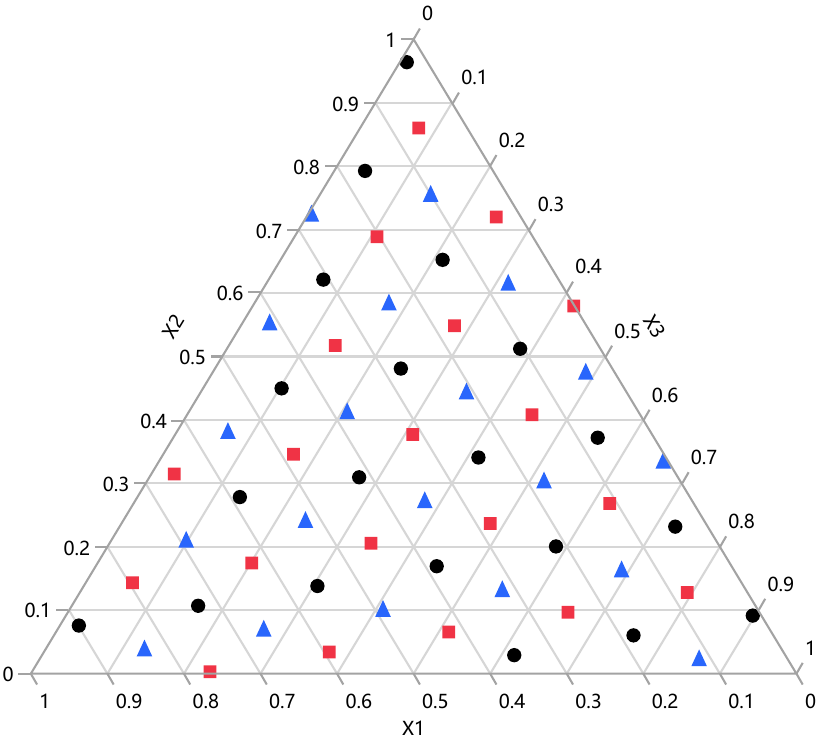}
        \caption{full design}
     \end{subfigure}
    \begin{subfigure}[b]{0.49\textwidth}
        \includegraphics[width=\textwidth]{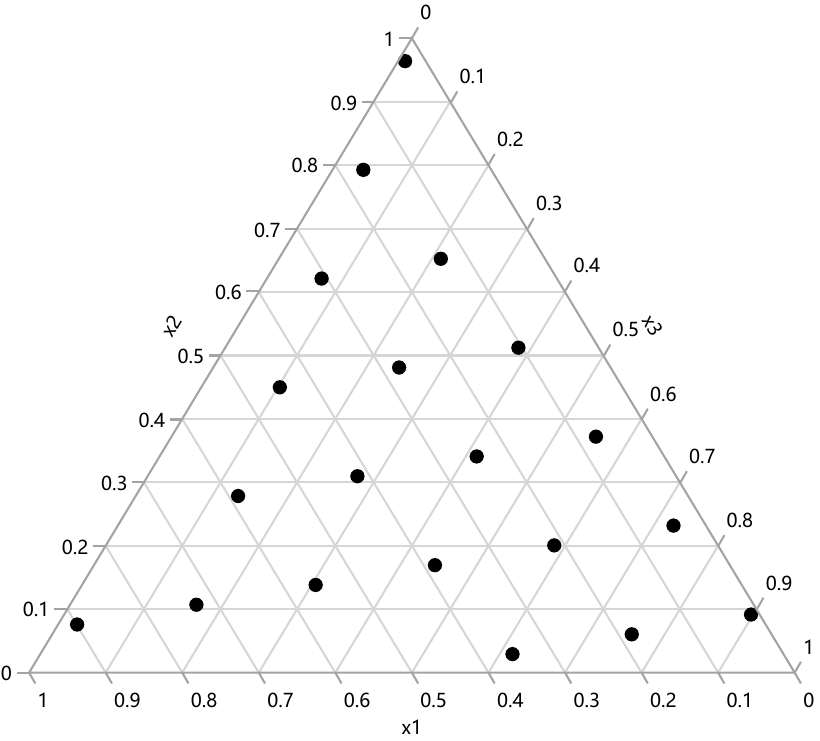}
        \caption{ slice 1}
  \end{subfigure}
  \begin{subfigure}[b]{0.49\textwidth}
        \includegraphics[width=\textwidth]{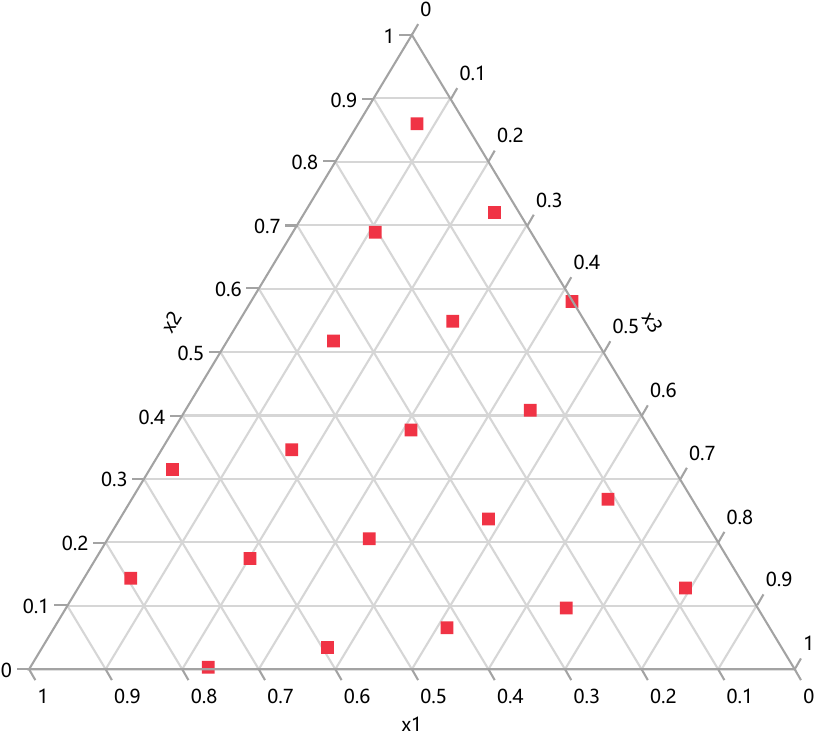}
        \caption{slice 2}
  \end{subfigure}
    \begin{subfigure}[b]{0.49\textwidth}
        \includegraphics[width=\textwidth]{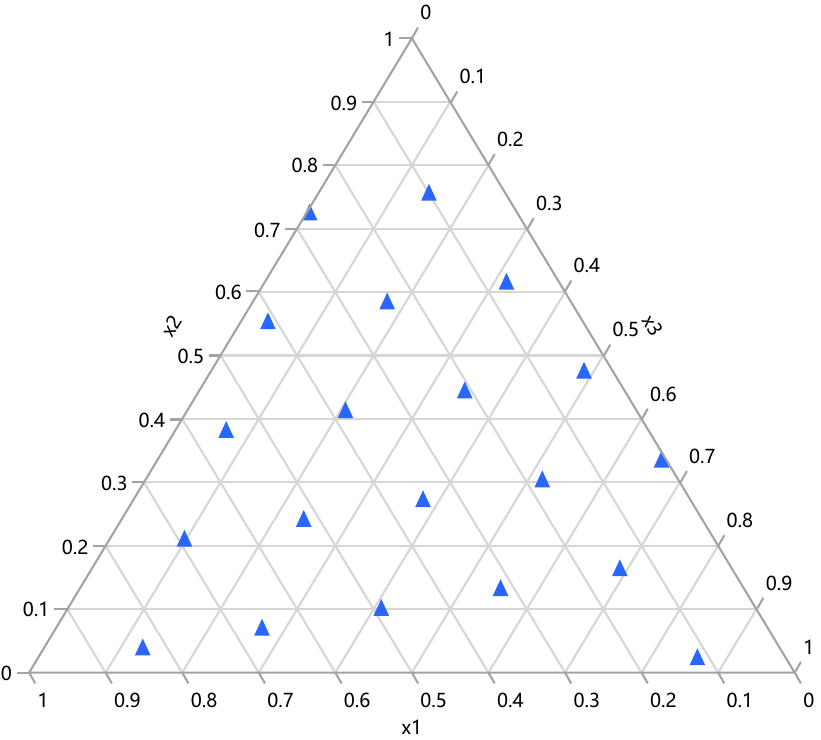}
        \caption{slice 3}
  \end{subfigure}
 \caption{\ Scatter plot of sliced space-filling design with mixtures generated by CDM}
\label{figCDM}
\end{figure}
\begin{figure}[htbp]
    \centering
    \begin{subfigure}[b]{0.49\textwidth}
        \includegraphics[width=\textwidth]{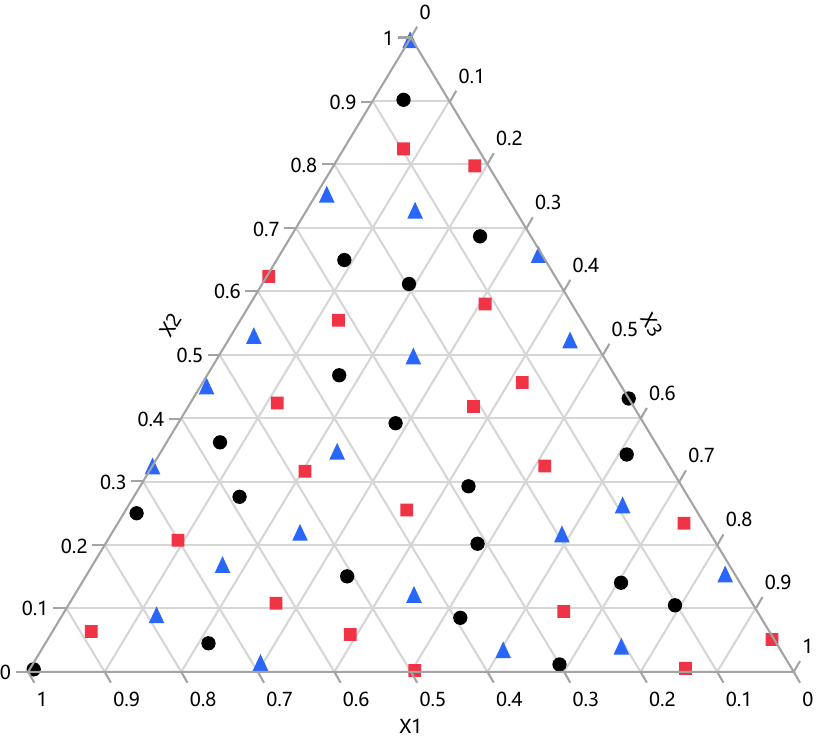}
        \caption{full design}
     \end{subfigure}
    \begin{subfigure}[b]{0.49\textwidth}
        \includegraphics[width=\textwidth]{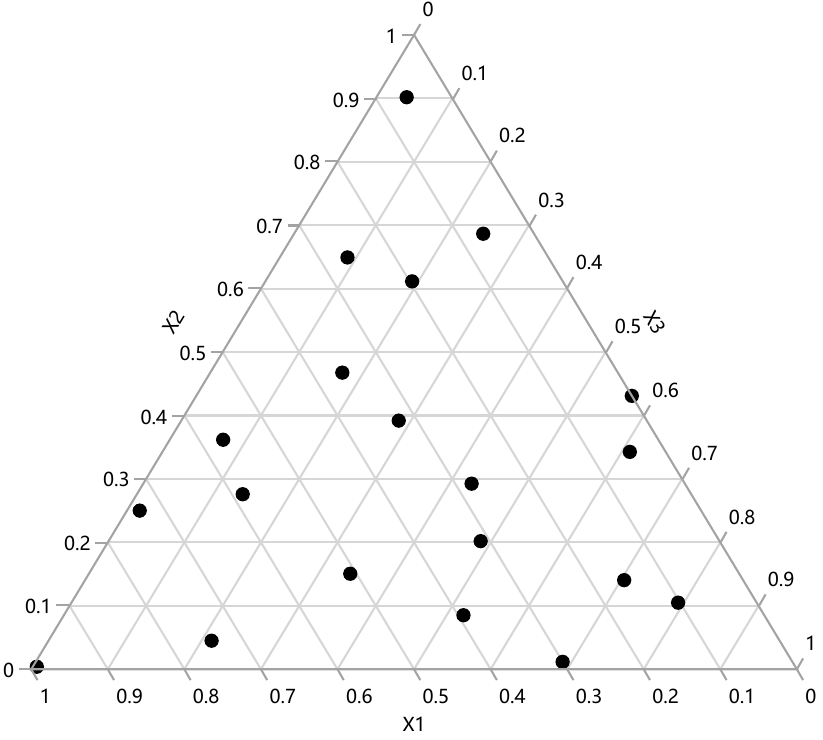}
        \caption{ slice 1}
  \end{subfigure}
  \begin{subfigure}[b]{0.49\textwidth}
        \includegraphics[width=\textwidth]{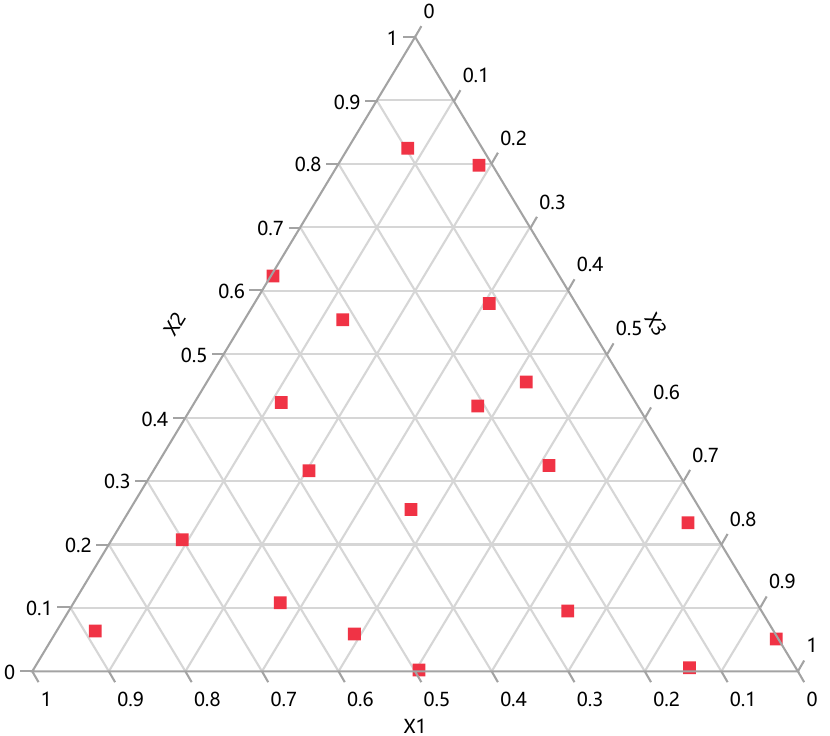}
        \caption{slice 2}
  \end{subfigure}
    \begin{subfigure}[b]{0.49\textwidth}
        \includegraphics[width=\textwidth]{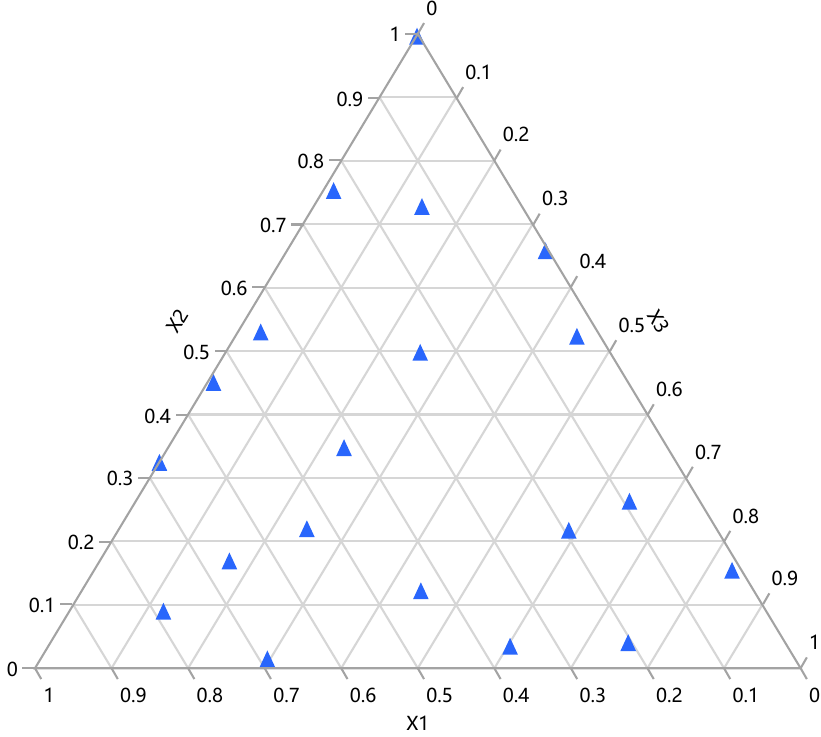}
        \caption{slice 3}
  \end{subfigure}
 \caption{\ Scatter plot of sliced space-filling design with mixtures generated by FFFM}
\label{figFFFM}
\end{figure}
\begin{figure}[htbp]
    \centering
    \begin{subfigure}[b]{0.49\textwidth}
        \includegraphics[width=\textwidth]{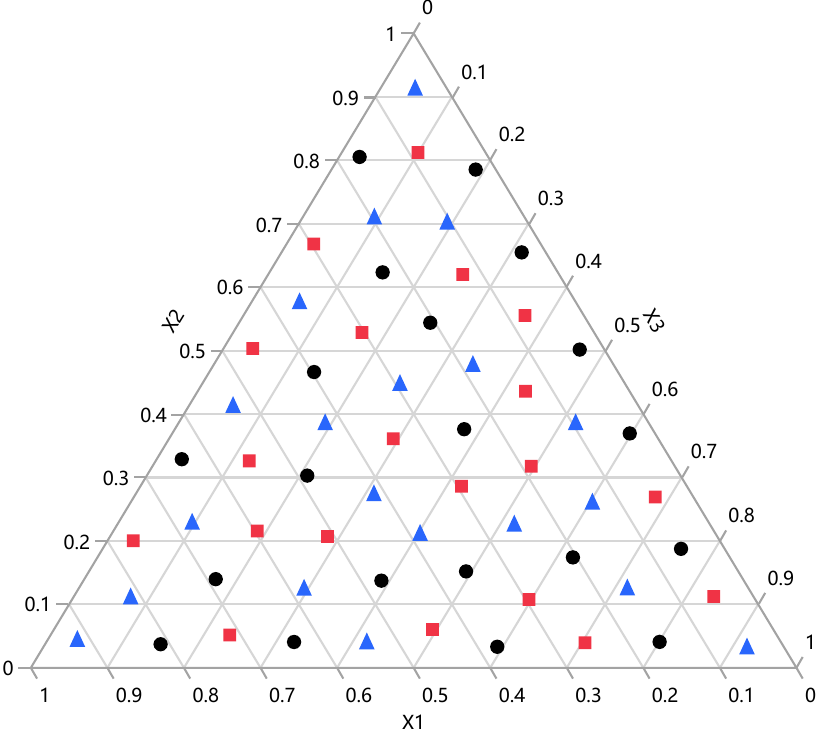}
        \caption{full design}
     \end{subfigure}
    \begin{subfigure}[b]{0.49\textwidth}
        \includegraphics[width=\textwidth]{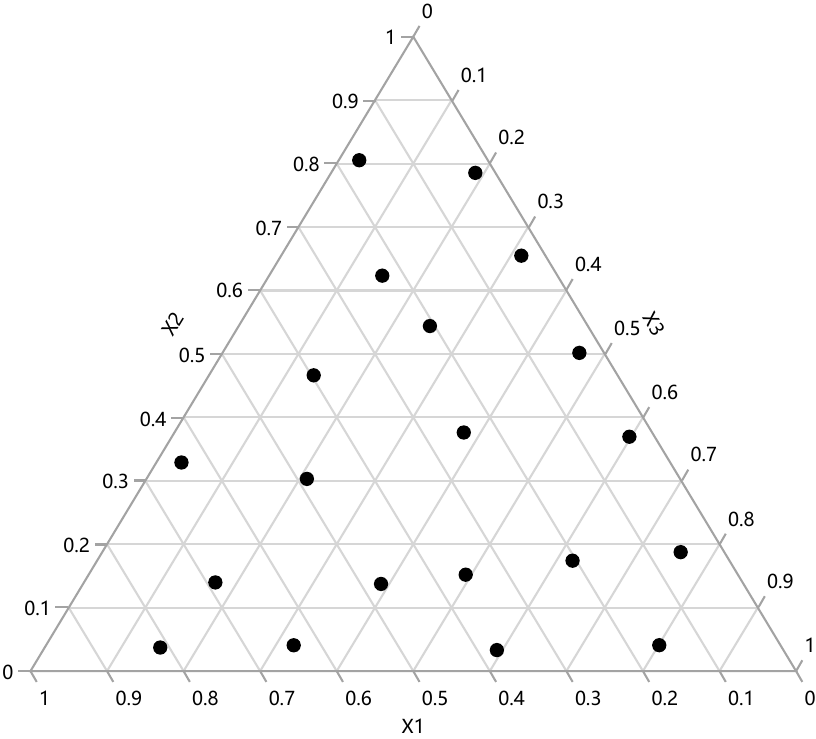}
        \caption{ slice 1}
  \end{subfigure}
  \begin{subfigure}[b]{0.49\textwidth}
        \includegraphics[width=\textwidth]{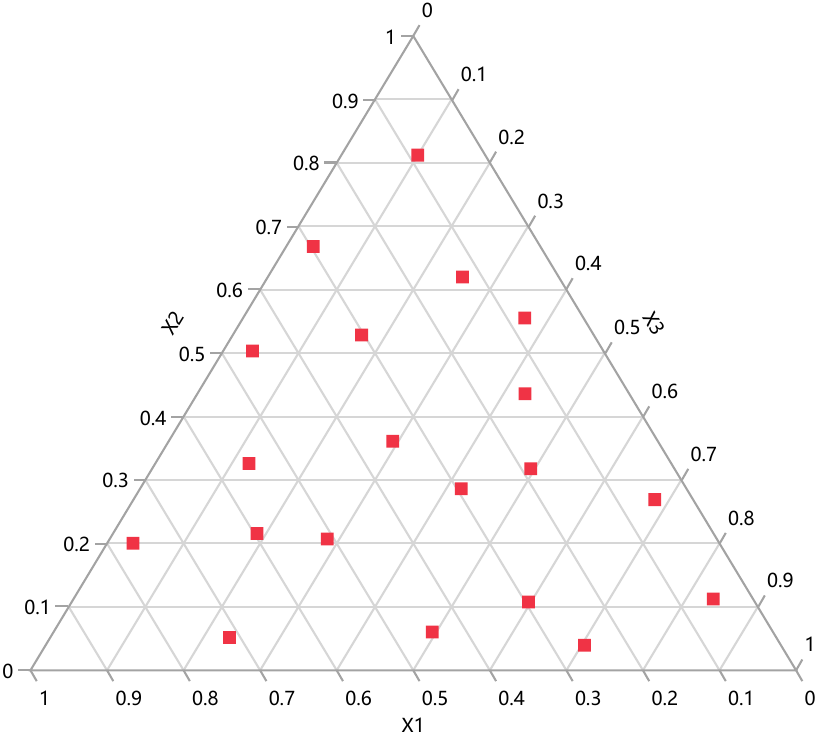}
        \caption{slice 2}
  \end{subfigure}
    \begin{subfigure}[b]{0.49\textwidth}
        \includegraphics[width=\textwidth]{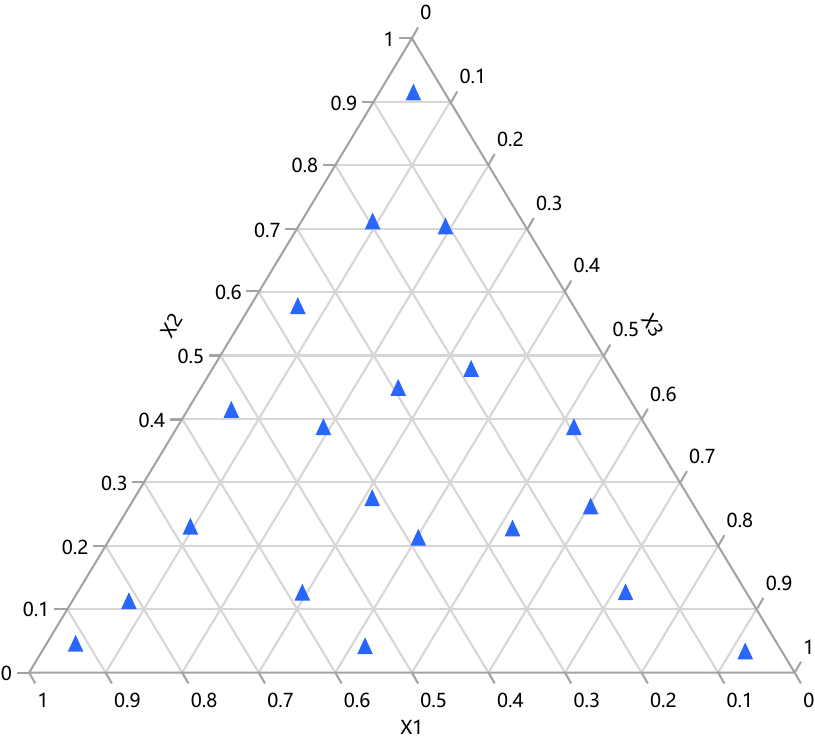}
        \caption{slice 3}
  \end{subfigure}
 \caption{\ Scatter plot of sliced space-filling design with mixtures generated by ParM}
\label{figParM}
\end{figure}
\begin{figure}[htbp]
    \centering
    \begin{subfigure}[b]{0.49\textwidth}
        \includegraphics[width=\textwidth]{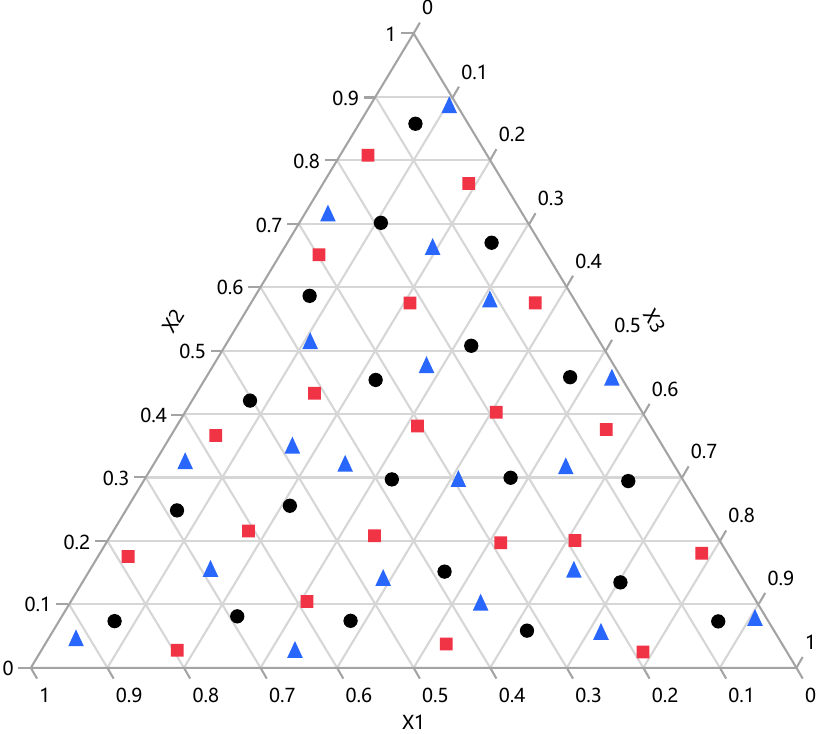}
        \caption{full design}
     \end{subfigure}
    \begin{subfigure}[b]{0.49\textwidth}
        \includegraphics[width=\textwidth]{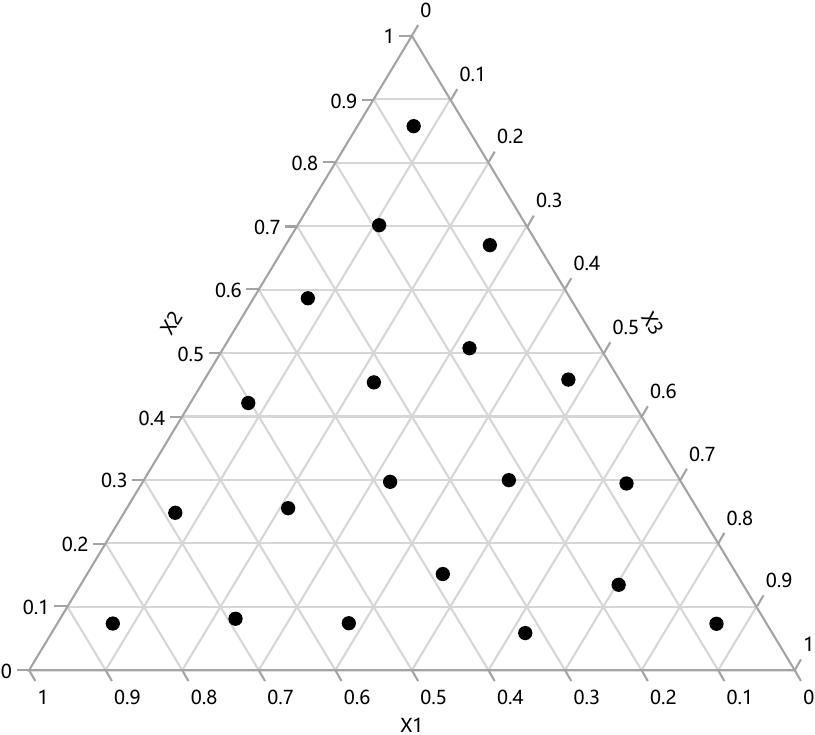}
        \caption{ slice 1}
  \end{subfigure}
  \begin{subfigure}[b]{0.49\textwidth}
        \includegraphics[width=\textwidth]{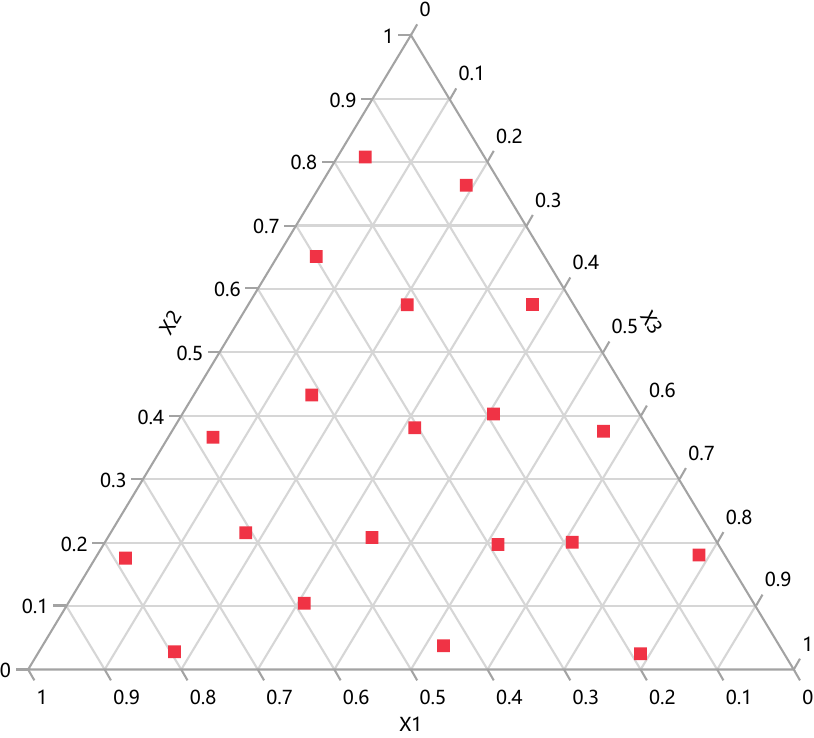}
        \caption{slice 2}
  \end{subfigure}
    \begin{subfigure}[b]{0.49\textwidth}
        \includegraphics[width=\textwidth]{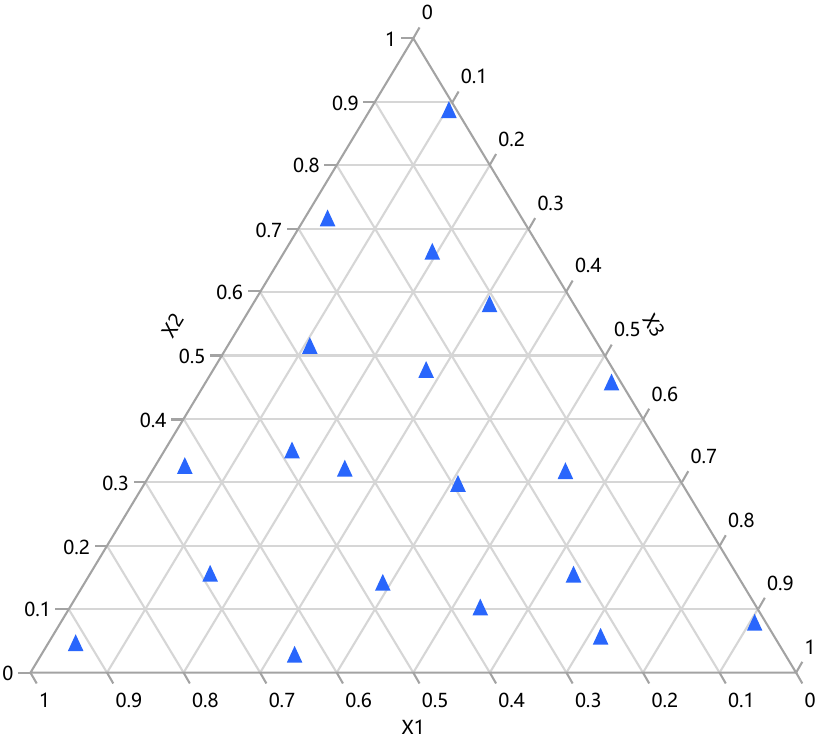}
        \caption{slice 3}
  \end{subfigure}
 \caption{\ Scatter plot of sliced space-filling design with mixtures generated by SeqM}
\label{figSeqM}
\end{figure}
\begin{figure}[htbp]
    \centering
    \begin{subfigure}[b]{0.49\textwidth}
        \includegraphics[width=\textwidth]{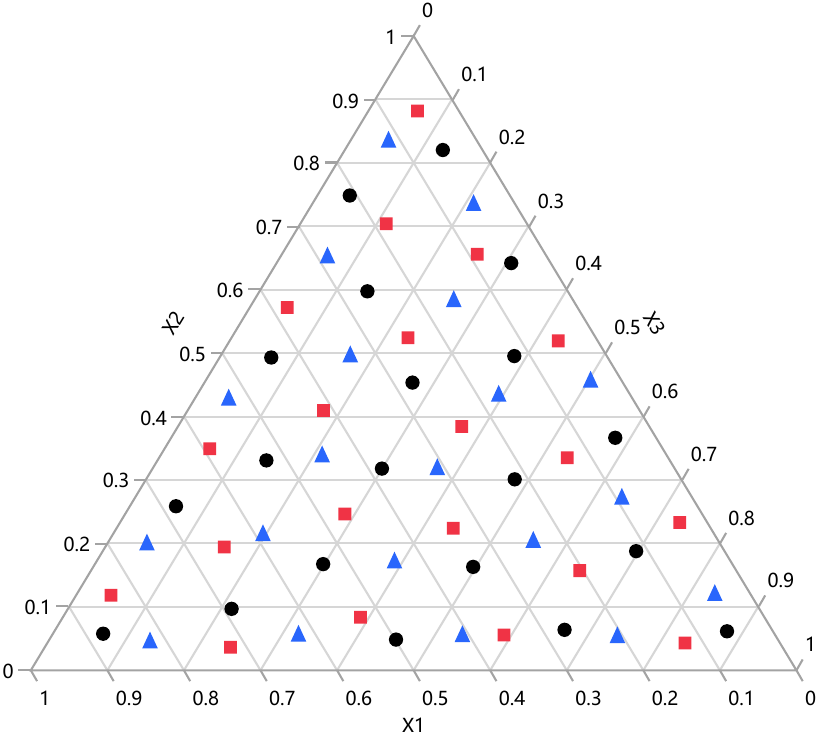}
        \caption{full design}
     \end{subfigure}
    \begin{subfigure}[b]{0.49\textwidth}
        \includegraphics[width=\textwidth]{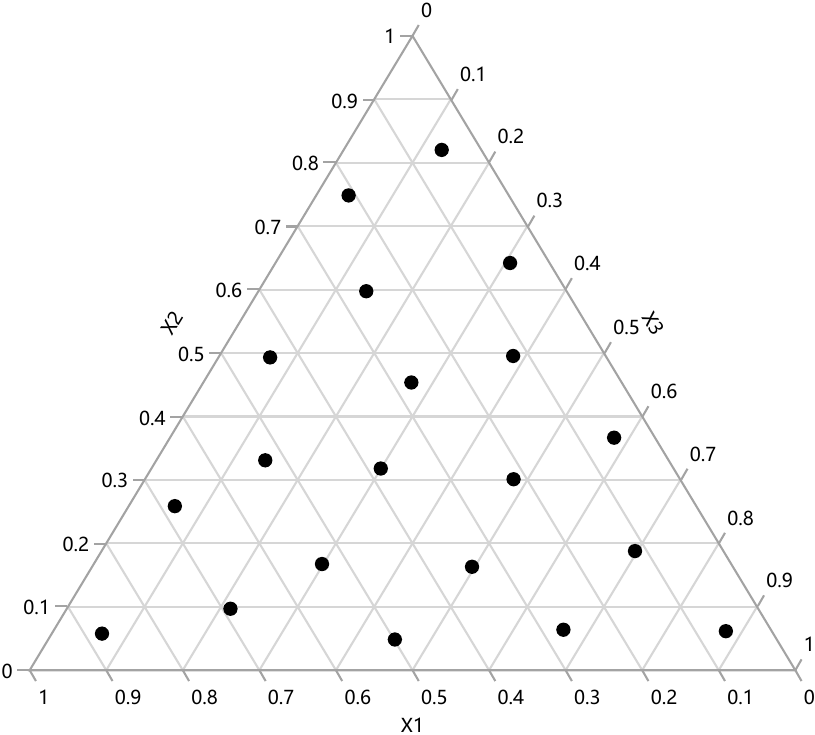}
        \caption{ slice 1}
  \end{subfigure}
  \begin{subfigure}[b]{0.49\textwidth}
        \includegraphics[width=\textwidth]{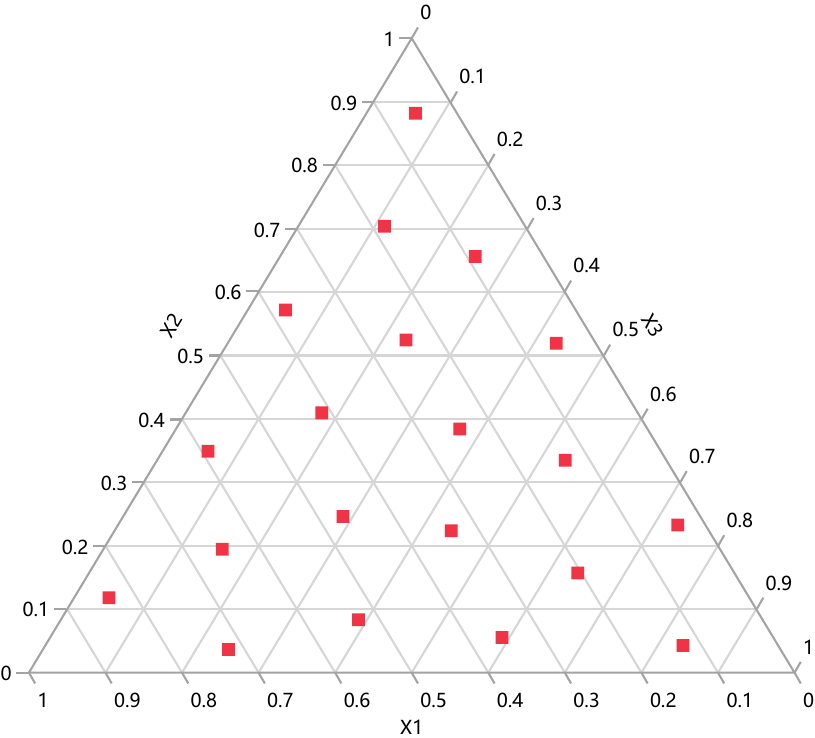}
        \caption{slice 2}
  \end{subfigure}
    \begin{subfigure}[b]{0.49\textwidth}
        \includegraphics[width=\textwidth]{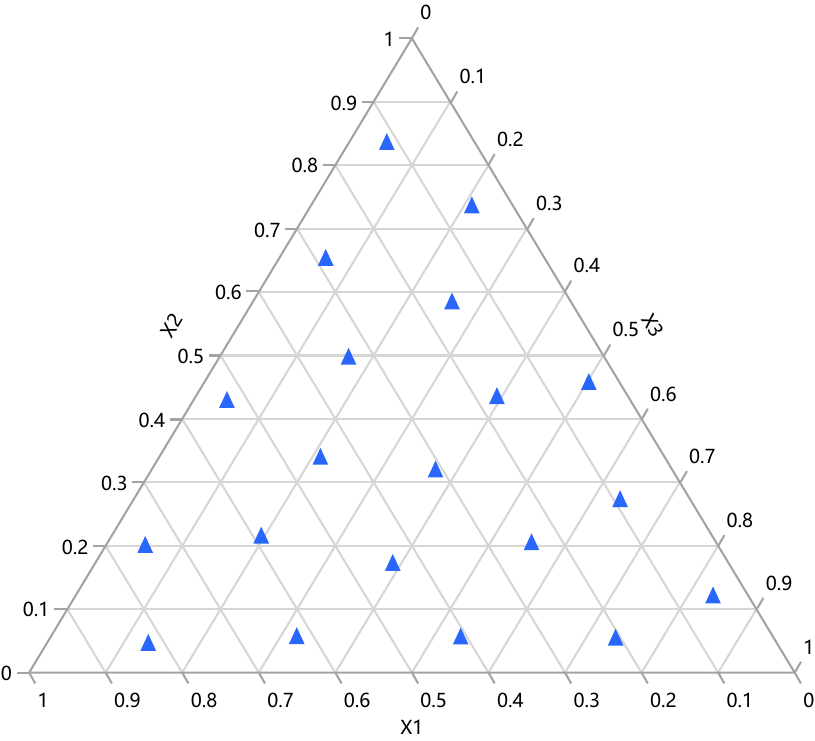}
        \caption{slice 3}
  \end{subfigure}
 \caption{\ Scatter plot of sliced space-filling design with mixtures generated by MHED}
\label{figMHED}
\end{figure}
\begin{figure}[htbp]
    \centering
    \begin{subfigure}[b]{0.49\textwidth}
        \includegraphics[width=\textwidth]{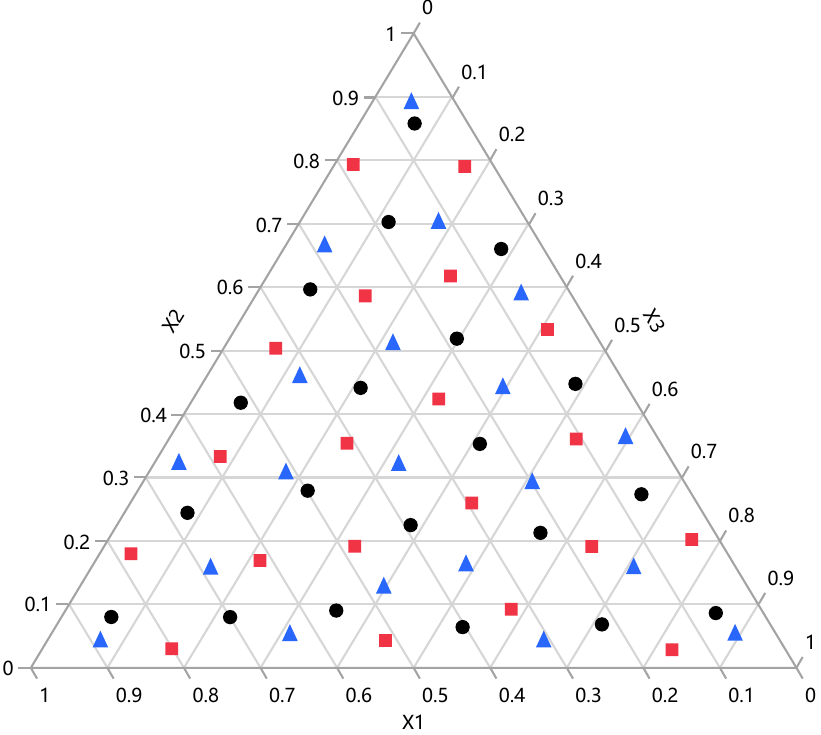}
        \caption{full design}
     \end{subfigure}
    \begin{subfigure}[b]{0.49\textwidth}
        \includegraphics[width=\textwidth]{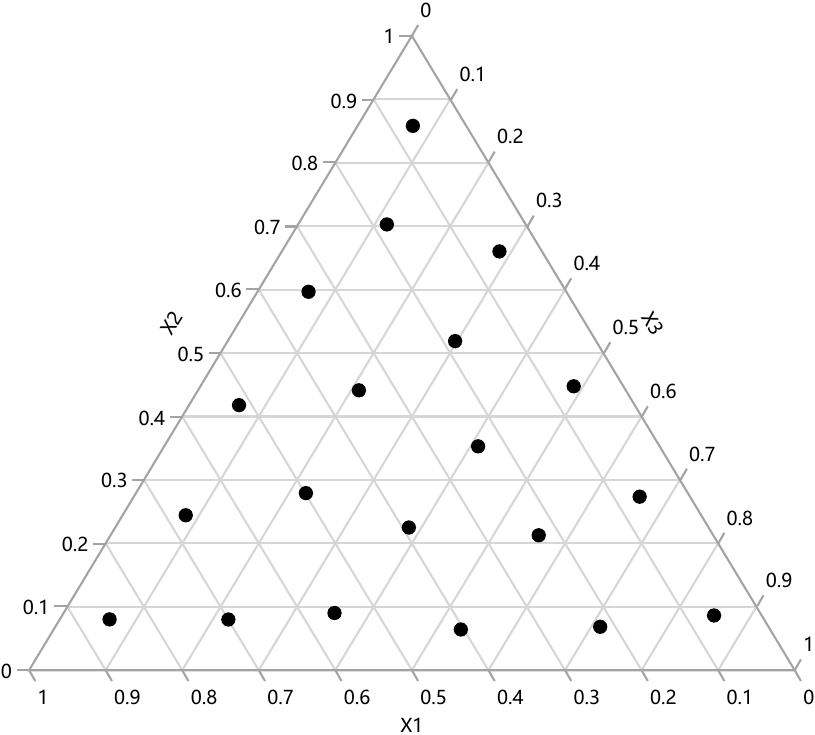}
        \caption{ slice 1}
  \end{subfigure}
  \begin{subfigure}[b]{0.49\textwidth}
        \includegraphics[width=\textwidth]{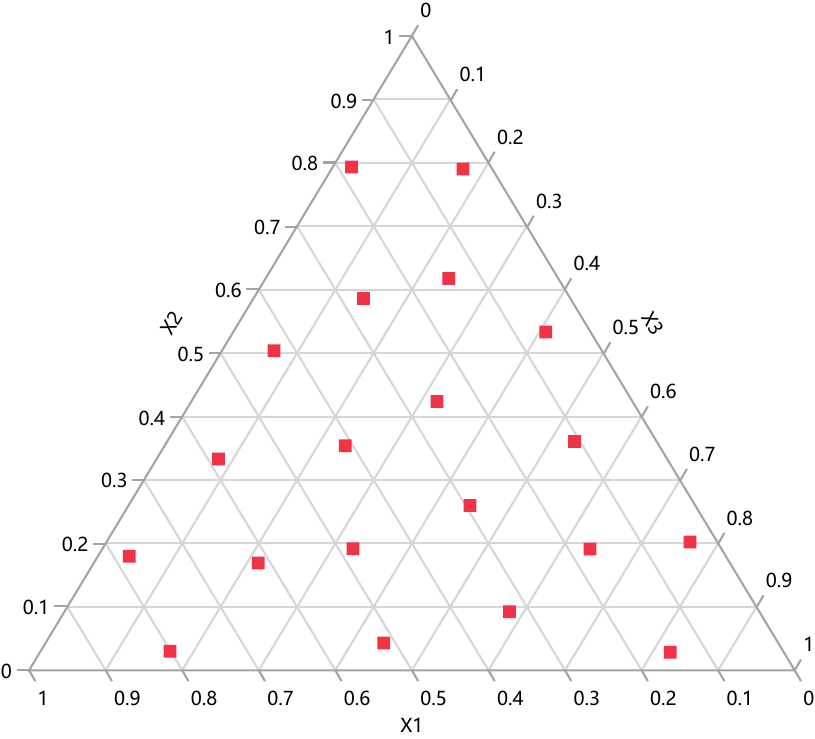}
        \caption{slice 2}
  \end{subfigure}
    \begin{subfigure}[b]{0.49\textwidth}
        \includegraphics[width=\textwidth]{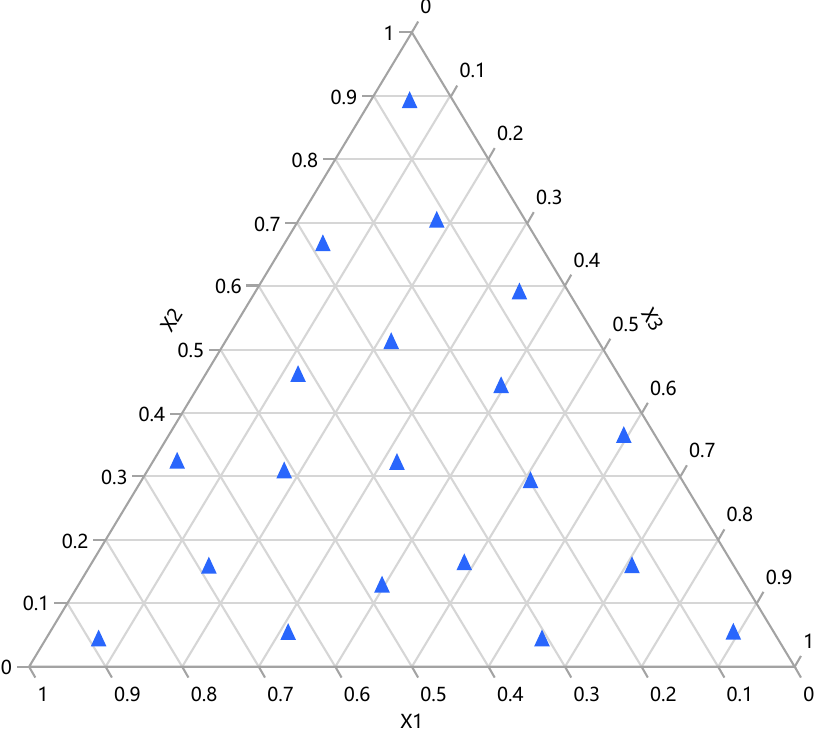}
        \caption{slice 3}
  \end{subfigure}
 \caption{\ Scatter plot of sliced space-filling design with mixtures generated by SeqHED}
\label{figSeqHED}
\end{figure}

\bibliography{bibliography.bib}